\documentclass[12pt]{article}
\usepackage{latexsym,amsfonts,amssymb}
\usepackage[dvips]{color}
\usepackage{colordvi,multicol}

\def \cal{\mathcal}

\newtheorem{thm}{Theorem}[section]

\newtheorem{lem}[thm]{Lemma}

\newtheorem{rem}[thm]{Remark}

\begin{document}
\title{\bf  Two Stronger Versions of  the  Union-closed Sets Conjecture }
\author {
Zhen Cui \\
{\small School of Mathematics and Statistics, Shandong Normal University, China }  \\ \\
Ze-Chun Hu\thanks{Corresponding author: College of Mathematics, Sichuan University, 610065, China\vskip 0cm E-mail addresses: cuizhen@sdnu.edu.cn (Z.
Cui), zchu@scu.edu.cn (Z.-C. Hu)}\\
  {\small  College of Mathematics, Sichuan University, China}}

\date{}
 \maketitle

\noindent{\bf Abstract}\quad  The union-closed sets conjecture (Frankl's conjecture)  says that  for any finite union-closed family of finite sets, other than the family
consisting only of the empty set, there exists an element that belongs to at least half of the sets in the family. In this paper, we introduce two stronger versions of Frankl's
 conjecture and give a partial proof. Three related questions are introduced.

\noindent{\bf Key words} The union-closed sets conjecture, Frankl's conjecture

\noindent{\bf Mathematics Subject Classification (2010)} 03E05, 05A05

\section{Introduction}

A family $\mathcal{A}$ of sets is union-closed if for every two member-sets $A,B\in \mathcal{A}$, their
union $A\cup  B$ is contained in $\mathcal{A}$. For simplicity, denote $n=|\cup_{A\in \cal{F}}A|$ and $m=|\cal{F}|$.

In 1979, Peter Frankl (cf. \cite{Ri85, St86}) conjectured that  for any finite union-closed family of finite sets, other than the family consisting only of the empty set, there exists an element that belongs to at least half of the sets in the family.

If a union-closed family $\cal{F}$ contains a set with one element or two elements, then Frankl's conjecture holds for $\cal{F}$ (\cite{SR89} ).
This result was extended by Poonen (\cite{Po92}). In addition,  the author in \cite{Po92}  proved that Frankl's conjecture holds if $n\leq 7$ or $m\leq 28$, and proved an equivalent lattice formulation of Frankl's conjecture. Bo\v{s}njak and Markovi\'{c} (\cite{BM08}) proved that Frankl's conjecture holds if $n\leq 11$.   Zivkovi\'{c} and Vu\v{c}kovi\'{c} (\cite{ZV12}) gave a computer
assisted proof that Frankl's conjecture is true if  $n\leq 12$, which together with Faro's result (\cite{Lo94-b}) (see also Roberts and Simposon \cite{RS10}) implies that Frankl's conjecture holds if  $m\leq 50$.  For more progress on Frankl's conjecture, we refer to \cite{BS15}, \cite{GY98}, \cite{JV98}, \cite{Lo94-a}, \cite{Ma07}, \cite{Mo06}, \cite{Va02}, \cite{Va03}, \cite{Va04}.

Let $M_n=\{1,2,\cdots,n\}$ and $\cal{F}\subset 2^{M_n}=\{A: A\subset M_n\}$ with $\cup_{A\in\cal{F}}A=M_n$.  Suppose that $\cal{F}$ is union-closed. We assume that $\emptyset \in \cal{F}$.  For any $k=1,2,\cdots,n$, denote $\cal{M}_k=\{A\in 2^{M_n}: |A|=k\}.$ Define
$$
T(\cal{F})=\inf\{1\leq k\leq n: \cal{F}\cap \cal{M}_k\neq \emptyset\}.
$$
Then $1\leq T(\cal{F})\leq n$.  By virtue of $T(\cal{F})$, we introduce the following two stronger versions of Frankl's conjecture.

{\bf $S_1$-Frankl conjecture:} If  $n\geq 2$ and $T(\cal{F})=k\in \{2,\cdots,n\}$, then there exist at least $k$ elements in $M_n$ which belong to at least half of the sets in $\cal{F}$.

{\bf $S_2$-Frankl conjecture:} If  $n\geq 2$ and $T(\cal{F})\geq 2$, then there exist at least two elements in $M_n$ which belong to at least half of the sets in $\cal{F}$.

\begin{rem}\label{rem1.1}  If $T(\cal{F})=1$, then by \cite{SR89}, we know that there exists an element that belongs to at least half of the sets in $\cal{F}$. Hence we have
\begin{eqnarray*}
S_1\mbox{-Frankl conjecture} \Rightarrow S_2\mbox{-Frankl conjecture} \Rightarrow \mbox{Frankl's conjecture}.
\end{eqnarray*}
\end{rem}

As to the case $T(\cal{F})=2$, we have the following lemma, which will be used in Section 3.

\begin{lem}\label{lem-1.2}
If $T(\cal{F})=2$ and there exist two sets $A_1$ and $A_2$ in $\cal{F}\cap \cal{M}_2$ such that $A_1\cap A_2=\emptyset$, then $S_1$-Frankl conjecture holds, i.e. there are at least two elements in $M_n$ which belong to at least half of the sets in $\cal{F}$.
\end{lem}
{\bf Proof.} By \cite{SR89}, we know that for $i=1,2$, there exists one element in $A_i$ that belongs to at least half of the sets in $\cal{F}$.
The condition that $A_1\cap A_2=\emptyset$ implies the claim of the lemma.\hfill\fbox

In next section, we will  give the proofs for $S_1$-Frankl conjecture for three cases: $T(\cal{F})=n, T(\cal{F})=n-1$ and $T(\cal{F})=n-2$. It follows that both $S_1$-Frankl conjecture and $S_2$-Frankl conjecture  hold if $2\leq n\leq 4$.  In Section 3, we will show that $S_1$-Frankl conjecture holds if $n=5$ and thus $S_2$-Frankl conjecture holds if $n=5$.  In the final section, a counterexample will be given for $S_1$-Frankl conjecture when $n=9$ (this example was introduced by an anonymous referee).


\section{Three cases on  $S_1$-Frankl conjecture}

In this section, we will  give the proofs for $S_1$-Frankl conjecture for three cases: $T(\cal{F})=n, T(\cal{F})=n-1$ and $T(\cal{F})=n-2$, respectively,  which will be used in next section.

\subsection{$T(\cal{F})=n\ (n\geq 2)$}

In this case,  $\cal{F}=\{\emptyset, M_n\}$. Then all the $n$ elements in $M_n$  belong to half of the sets in $\cal{F}$.
It follows that $S_1$-Frankl conjecture holds if $n=2$.

\subsection{$T(\cal{F})=n-1$  $(n\geq 3)$}

Now,  $\cal{F}=\{\emptyset, M_n\}\cup \cal{G}_{n-1}$, where $\cal{G}_{n-1}$ is a nonempty subset of $\cal{M}_{n-1}$.
 We will show that at least $n-1$ elements in $M_n$ belong to at least half of the sets in $\cal{F}$. We have two cases: $|\cal{G}_{n-1}|=1$ and $|\cal{G}_{n-1}|\geq 2$.

\begin{itemize}
\item[(1)]  $|\cal{G}_{n-1}|=1$.  Denote  $\cal{G}_{n-1}=\{G\}$. Then all the $n -1$ elements in $G$ belong to at least half of the sets in $\cal{F}$.

\item[(2)]   $|\cal{G}_{n-1}|\geq 2$. Note that for any element $i$  in $M_n$, it belongs to all the sets in $\cal{M}_{n-1}$ except one set $M_n\backslash \{i\}$.  Thus in this case, all the $n$ elements in $M_n$ belong to at least $|\cal{G}_{n-1}|-1$ set(s) in $\cal{G}_{n-1}$ and thus all the $n$ elements in $M_n$ belong to at least half of the sets in $\cal{G}_{n-1}$. It follows that all the $n$ elements in $M_n$ belong to at least half of the sets in $\cal{F}$.
\end{itemize}

By 2.1 and 2.2, we get that $S_1$-Frankl conjecture holds if $n=3$.

\subsection{$T(\cal{F})=n-2$  $(n\geq 4)$}

In the following, we will show  that at least $n-2$ elements in $M_n$ belong to at least half of the sets in $\cal{F}$.  We have two cases: $T(\cal{F}\backslash \cal{M}_{n-2})=n$ and $T(\cal{F}\backslash \cal{M}_{n-2})=n-1$.

\subsubsection{$T(\cal{F}\backslash \cal{M}_{n-2})=n$}

  Now $\cal{F}=\{\emptyset, M_n\}\cup \cal{G}_{n-2}$, where $\cal{G}_{n-2}$ is a nonempty subset of $\cal{M}_{n-2}$. We have two subcases: $|\cal{G}_{n-2}|=1$ and $|\cal{G}_{n-2}|\geq 2$.

\begin{itemize}
\item[(1)]  $|\cal{G}_{n-2}|=1$. Denote  $\cal{G}_{n-2}=\{G\}$. Then all the $n-2$ elements in $G$ belong to at least half of the sets in $\cal{F}$.

\item[(2)]   $|\cal{G}_{n-2}|=m\geq 2$.    Denote $\cal{G}_{n-2}=\{G_1,\cdots,G_m\}$.  Then for any $i,j=1,\cdots,m,i\neq j$, it holds that
$G_i\cup G_j=M_n$ and thus all the $n$ elements in $M_n$ belong to at least one of the two sets  $G_i$ and $G_j$. We have the following two subcases:

\begin{itemize}
\item[(2.1)]   $m$ is an even number. Then we have $G_1\cup G_2=\cdots=G_{m-1}\cup G_m.$ Hence  all the $n$ elements in $M_n$ belong to at least half of the sets in $\cal{G}_{n-2}$ and thus belong to at least half of the sets in $\cal{F}$.

\item[(2.2)] $m$ is an odd number. Now for any $i=1,\cdots,m$, we know that $|\cal{G}_2\backslash \{G_i\}|$ is an even number. By (2.1), we know that all the $n$ elements in $M_n$ belong to at least half of the sets in $\cal{G}_2\backslash \{G_i\}$ and thus all the $n-2$ elements in $G_i$ belong to at least half of the sets in $\cal{F}$. Since $G_1\cup G_2=M_n$, we obtain that all the $n$ elements in $M_n$ belong to at least half of the sets in $\cal{F}$.

 \end{itemize}
\end{itemize}

\subsubsection{$T(\cal{F}\backslash \cal{M}_{n-2})=n-1$}

  Now $\cal{F}=\{\emptyset, M_n\}\cup \cal{G}_{n-2}\cup \cal{G}_{n-1}$, where $\cal{G}_{n-i}$ is a nonempty subset of $\cal{M}_{n-i}$ for $i=1,2$.  We have two subcases: $|\cal{G}_{n-1}|=1$ and $|\cal{G}_{n-1}|\geq 2$.

\begin{itemize}
\item[(1)]  $|\cal{G}_{n-1}|=1$. Without loss of generality, we assume that  $\cal{G}_{n-1}=\{\{1,2,\cdots,n-1\}\}$.  We have two subcases: $|\cal{G}_{n-2}|=1$ and $|\cal{G}_{n-2}|\geq 2$.

\begin{itemize}
\item[(1.1)] $|\cal{G}_{n-2}|=1$. Denote $\cal{G}_{n-2}=\{G\}$. Then all the elements in $G\cup \{1,2,\cdots,n-1\}$ belong to at least one of the two sets in $\cal{G}_{n-2}\cup \cal{G}_{n-1}$ and thus belong to at  least  half of the sets in $\cal{F}$.

\item[(1.2)] $|\cal{G}_{n-2}|=m\geq 2$. Denote $\cal{G}_{n-2}=\{G_1,\cdots,G_m\}$. Then for any $i,j=1,\cdots,m,i\neq j$, either  $G_i\cup G_j=M_n$ or $G_i\cup G_j=\{1,2,\cdots,n-1\}$.
Now we have the following three subcases:

\begin{itemize}
\item[(1.2.1)] $m$ is an even number and there exists a permutation $(i_1,\cdots,i_m)$ of $(1,2,\cdots,m)$ such that
$
G_{i_1}\cup G_{i_2}=\cdots=G_{i_{m-1}}\cup G_{i_m}=M_n.
$
Then all the $n$ elements in $M_n$ belong to at least half of the sets in $\cal{G}_{n-2}$ and thus all the $n-1$ elements in $\{1,2,\cdots,n-1\}$ belong to at least half of the sets in $\cal{F}$.

\item[(1.2.2)] $m$ is an odd number and there exists a permutation $(i_1,\cdots,i_m)$ of $(1,2,\cdots,m)$ such that
$
G_{i_1}\cup G_{i_2}=\cdots=G_{i_{m-2}}\cup G_{i_{m-1}}=M_n.
$
Then all the $n$ elements in $M_n$ belong to at least half of the sets in $\{G_{i_1},\cdots,G_{i_{m-1}}\}$.  Note that all the elements in $G_{i_m}\cup \{1,2,\cdots,n-1\}$ belong to at least one of the two sets $G_{i_m}$ and $\{1,2,\cdots,n-1\}$. Then we get that all the elements in $G_{i_m}\cup \{1,2,\cdots,n-1\}$ belong to at least half of the sets in $\cal{F}$.

\item[(1.2.3)] We can decompose $\cal{G}$ into two disjoint parts $\{G_{i_1},\cdots,G_{i_{2k}}\}$ (hereafter, this part may be an empty set) and $\{G_{i_{2k+1}},\cdots,G_{i_m}\}$, where
$\{i_1,\cdots,i_m\}=\{1,2,\cdots,m\}$, $m-2k\geq 2$, and

\quad (i) $G_{i_1}\cup G_{i_2}=\cdots=G_{i_{2k-1}}\cup G_{i_{2k}}=M_n$;

\quad (ii)  for any two different indexes $i,j$ from $\{i_{2k+1},\cdots,i_m\}$, $G_i\cup G_j=\{1,2,\cdots,n-1\}$.\\
Then all the $n$ elements in $M_n$ belong to at least half of the sets in $\{G_{i_1},\cdots,G_{i_{2k}}\}$. As to the second part $\{G_{i_{2k+1}},\cdots,G_{i_m}\}$, by following Case (2) in Section 2.2, we know that all the $n-1$ elements in $\{1,2,\cdots,n-1\}$  belong to at least\linebreak  $|\{G_{i_{2k+1}},\cdots,G_{i_m}\}|-1$ set(s) in $\{G_{i_{2k+1}},\cdots,G_{i_m}\}$. Thus all the $n-1$ elements in $\{1,2,\cdots, n-1\}$ belong to at least half of the sets in $\{G_{i_{2k+1}},\cdots,G_{i_m}\}$.  It follows that all the $n-1$ elements in $\{1,2,\cdots, n-1\}$ belong to at least half of the sets in $\cal{F}$.
\end{itemize}
\end{itemize}

\item[(2)]  $|\cal{G}_{n-1}|\geq 2$. Now by Case (2) in Section 2.2, we know that  all the $n$ elements in $M_n$ belong to at least  $|\cal{G}_{n-1}|-1$ set(s) in $\cal{G}_{n-1}$.  It follows that all the $n$ elements in $M_n$ belong to at least  half of the sets in $\cal{G}_{n-1}$ and thus it is enough to show that at least $n-2$ elements in $M_n$ belong to at least half of the sets in $\cal{G}_{n-2}$ or $\cal{G}_{n-2}\cup \cal{G}_{n-1}$. We have two subcases: $|\cal{G}_{n-2}|=1$ and $|\cal{G}_{n-2}|\geq 2$.

\begin{itemize}
\item[(2.1)] $|\cal{G}_{n-2}|=1$. Denote  $\cal{G}_{n-2}=\{G\}$. Now all the $n-2$ elements in $G$ belong to the unique set in
 $\cal{G}_{n-2}$.

\item[(2.2)] $|\cal{G}_{n-2}|=m\geq 2$. Denote $\cal{G}_{n-2}=\{G_1,G_2,\cdots,G_m\}$.  For any $i,j=1,\cdots,m$, either $G_i\cup G_j=M_n$ or $|G_i\cup G_j|=n-1$. If  $G_i\cup G_j=M_n$, then all the $n$ elements in $M_n$ belong to at least one of the two sets $G_i$ and $G_j$.  We have the following three subcases:

\begin{itemize}
\item[(2.2.1)] $m$ is an even number and there exists a permutation $(i_1,\cdots,i_m)$ of $(1,2,\cdots,m)$ such that
$
G_{i_1}\cup G_{i_2}=\cdots=G_{i_{m-1}}\cup G_{i_m}=M_n.
$
Then all the $n$ elements in $M_n$ belong to at least half of the sets in $\cal{G}_{n-2}$.

\item[(2.2.2)] $m$ is an odd number and there exists a permutation $(i_1,\cdots,i_m)$ of $(1,2,\cdots,m)$ such that
$
G_{i_1}\cup G_{i_2}=\cdots=G_{i_{m-2}}\cup G_{i_{m-1}}=M_n.
$
Then all the $n$ elements in $M_n$ belong to at least half of the sets in $\{G_{i_1},\cdots,G_{i_{m-1}}\}$. Thus all the $n-2$ elements in  $G_{i_m}$ belong to at least half of the sets in $\cal{G}_{n-2}$.

\item[(2.2.3)]   We can decompose $\cal{G}_{n-2}$ into two disjoint parts $\{G_{i_1},\cdots,G_{i_{2k}}\}$ and $\{G_{i_{2k+1}},\cdots,G_{i_m}\}$, where $\{i_1,i_2,\cdots,i_m\}=\{1,2,\cdots,m\}, m-2k\geq 2$, and

  \quad (i) $G_{i_1}\cup G_{i_2}=\cdots=G_{i_{2k-1}}\cup G_{i_{2k}}=M_n$;

  \quad (ii) for any two different indexes $i,j$ from $\{i_{2k+1},\cdots,i_m\}$, $|G_i\cup G_j|=n-1$.\\
  Without loss of generality, we assume $G_{i_{2k+1}}=\{1,2,\cdots,n-2\}$.  Then we have the following three cases:

 \quad  (a)  $\cup_{j=2k+1}^mG_{i_j}=\{1,2,\cdots,n-2, n-1\}$;

  \quad (b)   $\cup_{j=2k+1}^mG_{i_j}=\{1,2,\cdots,n-2, n\}$;

 \quad  (c)  $\cup_{j=2k+1}^mG_{i_j}=\{1,2,\cdots,n-2, n-1, n\}$.
 \smallskip

  \quad As to Case (a), by  following Case (2) in Section 2.2, we know that all the $n-1$ elements in $\{1,2,\cdots,n-2,n-1\}$ belong to at least $|\{G_{i_{2k+1}},\cdots,G_{i_m}\}|-1$ set(s) in $\{G_{i_{2k+1}},\cdots,G_{i_m}\}$.  Hence in this case, all the $n-1$ elements in $\{1,2,\cdots,n-2,n-1\}$ belong to at least half of the sets in $\cal{G}_{n-2}$.

 \smallskip

  \quad Similarly, as to Case (b),  all the $n-1$ elements in $\{1,2,\cdots,n-2,n\}$ belong to at least half of the sets in $\cal{G}_{n-2}$.

 \smallskip

  \quad As to Case (c), define
  \begin{eqnarray*}
  &&\cal{H}_{n-1}:=\{A\in \{G_{i_{2k+2}},\cdots,G_{i_m}\}: A\cup G_{i_{2k+1}}=\{1,2,\cdots,n-2,n-1\}\},\\
  &&\cal{H}_{n}:=\{A\in \{G_{i_{2k+2}},\cdots,G_{i_m}\}: A\cup G_{i_{2k+1}}=\{1,2,\cdots,n-2,n\}\}.
  \end{eqnarray*}
 Suppose that $A\in \cal{H}_{n-1}$. Without loss of generality, we assume that $A=\{2,\cdots,n-2,n-1\}$.
Then by the above condition (ii) we know that $\cal{H}_{n}$ has a unique element $B=\{2,\cdots,n-2,n\}$, and further by the above condition (ii)  again we know that $A$ is the unique element of $\cal{H}_{n-1}$.  Thus in this case $m-2k=3$ and
$\{G_{i_{2k+1}},\cdots,G_{i_m}\}$ equals the following set
$$
\{\{1,2,\cdots,n-2\}, \{2,\cdots,n-2,n-1\}, \{2,\cdots,n-2,n\}\}.
$$
It follows that $\cal{G}_{n-1}$ at least contains the following three sets:
$$
\{1,2,\cdots,n-2,n-1\}, \ \{1,2,\cdots,n-2,n\}, \{2,\cdots,n-2, n-1,n\}.
$$
Hence all the $n-3$ elements in $\{2,\cdots,n-2\}$ belong to at least half of the sets in $\cal{G}_{n-2}$. As to the other three elements in $\{1,n-1,n\}$, they  belong to one of the three sets in $\{G_{i_{2k+1}},\cdots,G_{i_m}\}$ and belong to at least $|\cal{G}_{n-1}|-1$ set(s) in $\cal{G}_{n-1}$. It follows that all the three elements in $\{1,n-1,n\}$  belong to at least  $|\cal{G}_{n-1}|$ sets in $\cal{G}_{n-1}\cup \{G_{i_{2k+1}},\cdots,G_{i_m}\}$. Since $|\cal{G}_{n-1}|\geq 3$ and $ |\{G_{i_{2k+1}},\cdots,G_{i_m}|=3$, we obtain that  all the three elements in $\{1,n-1,n\}$  belong to at least half of the  sets in $\cal{G}_{n-1}\cup \{G_{i_{2k+1}},\cdots,G_{i_m}\}$. Hence in this case all the $n$ elements in $M_n$ belong to at least half of the sets in $\cal{G}_{n-1}\cup \cal{G}_{n-2}$ and thus belong to at least half of the sets in $\cal{F}$.
\end{itemize}
\end{itemize}
\end{itemize}

By 2.1, 2.2 and 2.3, we get that $S_1$-Frankl conjecture holds if $n=4$.

\section{$S_1$-Frankl conjecture for $n=5$}

In this section, we will prove that $S_1$-Frankl conjecture holds if $n=5$.   Let $M_5=\{1,2,\cdots,5\}$ and $\cal{F}\subset 2^{M_5}=\{A: A\subset M_5\}$ with $\cup_{A\in\cal{F}}A=M_5$.  Suppose that $\cal{F}$ is union-closed and  $\emptyset \in \cal{F}$.  For $k=1,2,\cdots,5$, denote $\cal{M}_k=\{A\in 2^{M_5}: |A|=k\}. $ Define
$
T(\cal{F})=\inf\{1\leq k\leq 5: \cal{F}\cap \cal{M}_k\neq \emptyset\}.
$
Then $1\leq T(\cal{F})\leq 5$.

By Section 2, we know that if $T(\cal{F})\in \{3,4, 5\}$, then there exist at least $T(\cal{F})$ elements in $M_5$ which belong to at least half of the sets in $\cal{F}$. Thus we need only consider the case $T(\cal{F})=2$ and will prove that there exist at least 2 elements which belong to at least half of the sets in $\cal{F}$.  We have three  subcases: $T(\cal{F}\backslash \cal{M}_2)=5,4,3$.

\subsection{$T(\cal{F}\backslash \cal{M}_2)=5$}

Now $\cal{F}=\{\emptyset, M_5\}\cup \cal{G}_2$, where $\cal{G}_2$ is a nonempty subset of $\cal{M}_2$.  If $A,B\in \cal{M}_2,A\neq B$, then $|A\cup B|\in \{3, 4\}$. It follows that $|\cal{G}_2|=1$.  Denote $\cal{G}_2=\{G\}$. Then all the 2 elements in $G$ belong to at least half of the sets in $\cal{F}$.

\subsection{$T(\cal{F}\backslash \cal{M}_2)=4$}

Now $\cal{F}=\{\emptyset, M_5\}\cup \cal{G}_2\cup \cal{G}_4$, where $\cal{G}_i$ is a nonempty subset of $\cal{M}_i$ for $i=2,4$.
We have two subcases: $|\cal{G}_2|\geq 2$ and $|\cal{G}_2|=1$.

\begin{itemize}
\item[(1)] $|\cal{G}_2|=m\geq 2$.  Denote $\cal{G}_2=\{G_1,\cdots,G_m\}$. Then for any $i,j=1,\cdots,m,i\neq j$, we must have that $G_i\cup G_j\in \cal{G}_4$,  and so  $G_i\cap G_j=\emptyset$.   By Lemma \ref{lem-1.2} we get that there are at least two elements in $M_5$ which belong to at least half of the sets in $\cal{F}$.

\item[(2)] $|\cal{G}_2|=1$. Denote $\cal{G}_2=\{G\}$. We have two subcases: $|\cal{G}_4|=1$ and $|\cal{G}_4|\geq 2$.

\begin{itemize}
\item[(2.1)]  $|\cal{G}_4|=1$. Denote $\cal{G}_4=\{H\}$. Now all the elements in $G\cup H$ belong to at least half of the sets in $\cal{F}$.

\item[(2.2)] $|\cal{G}_4|\geq 2$.  By Case (2) in Section 2.2, we know that all the 5 elements in $M_5$ belong to at least $|\cal{G}_4|-1$ set(s) in $\cal{G}_4$ and thus belong to at least half of the sets in $\cal{G}_4$. Hence all the 2 elements in $G$ belong to at least half of the sets in $\cal{F}$.
\end{itemize}
\end{itemize}

\subsection{$T(\cal{F}\backslash \cal{M}_2)=3$}

Now we have two subcases: $T(\cal{F}\backslash (\cal{M}_2\cup \cal{M}_3))=5$ and  $T(\cal{F}\backslash (\cal{M}_2\cup \cal{M}_3))=4$.

\subsubsection{$T(\cal{F}\backslash (\cal{M}_2\cup \cal{M}_3))=5$ }

Now $\cal{F}=\{\emptyset, M_5\}\cup \cal{G}_2\cup \cal{G}_3$, where $\cal{G}_i$ is a nonempty subset of $\cal{M}_i$ for $i=2,3$.
 We have two subcases: $|\cal{G}_3|=1$ and $|\cal{G}_3|\geq 2$.

\begin{itemize}
\item[(1)] $|\cal{G}_3|=1$. Without loss of generality, we assume $\cal{G}_3=\{\{1,2,3\}\}$. We have two subcases: $|\cal{G}_2|=1$ and $|\cal{G}_2|\geq 2$.

\begin{itemize}
\item[(1.1)] $|\cal{G}_2|=1$. Denote $\cal{G}_2=\{G\}$. Then all the elements in $G\cup \{1,2,3\}$ belong to at least half of the sets in $\cal{F}$.

\item[(1.2)] $|\cal{G}_2|\geq 2$. Denote $\cal{G}_2=\{G_1,\cdots,G_m\}$. Then for any $i,j=1,\cdots,m,i\neq j$, we must have that $G_i\cup G_j=\{1,2,3\}$. By Case (2) in Section 2.2, we know that all the 3 elements in $\{1,2,3\}$ belong to at least $|\cal{G}_2|-1$ set(s) in $\cal{G}_2$ and thus belong to at least half of the sets in $\cal{G}_2$. Hence all the 3 elements in $\{1,2,3\}$ belong to at least half of the sets in $\cal{F}$.
\end{itemize}

\item[(2)] $|\cal{G}_3|=m\geq 2$. Denote $\cal{G}_3=\{G_1,\cdots,G_m\}$. Then for any $i,j=1,\cdots,m,i\neq j$, we must have that
$G_i\cup G_j=M_5$ and thus all the 5 elements in $M_5$ belong to at least one of the two sets $G_i$ and $G_j$. We have the following two subcases:
\begin{itemize}
\item[(2.1)] $m$ is an even number. Now we have $G_1\cup G_2=\cdots=G_{m-1}\cup G_m=M_5$ and thus all the 5 elements in $M_5$ belong to at least half of the sets in $\cal{G}_3$. Hence it is enough to show that there exist 2 elements in $M_5$ which belong to at least half of the sets in $\cal{G}_2$. We have two subcases: $|\cal{G}_2|=1$ and $|\cal{G}_2|\geq 2$.

\begin{itemize}
\item[(2.1.1)] $|\cal{G}_2|=1$. Denote $\cal{G}_2=\{H\}$. Then the 2 elements in $H$ satisfy the condition.

\item[(2.1.2)] $|\cal{G}_2|=n\geq 2$. Denote $\cal{G}_2=\{H_1,\cdots,H_n\}$. Then for any $i,j=1,\cdots,n,i\neq j$, we must have that $H_i\cup H_j\in \cal{M}_3$. It follows that $n\leq 4$. Without loss of generality, we assume that $H_1=\{1,2\}$.

\begin{itemize}
\item[(2.1.2.1)] $n= 2$. Now $H_2\in \{\{1,3\}, \{1,4\}, \{1,5\}, \{2,3\}, \{2,4\}, \{2,5\}\}$.  Take $H_2=\{1,3\}$ for example. Then all the 3 elements in $H_1\cup H_2(=\{1,2,3\})$ belong to at least half of the sets in $\cal{G}_2$.

\item[(2.1.2.2)] $n= 3$. Now we have two subcases: $H_1\cap H_2\cap H_3=\emptyset$ and $|H_1\cap H_2\cap H_3|=1$.

\quad (a) $H_1\cap H_2\cap H_3=\emptyset$. Now $\cal{G}_2=\{\{1,2\},\{1,j\},\{2,j\}\}$ for $j=3,4,5$.  Take $j=3$ for example. Then we know that all the 3 elements in $\{1,2,3\}$ belong to at least two sets among the three sets in $\cal{G}_2$.
\smallskip

\quad (b) $|H_1\cap H_2\cap H_3|=1$.  Take $\cal{G}_2=\{\{1,2\}, \{1,3\}, \{1,4\}\}$ for example. Then $\{1,2,3,4\}\in \cal{F}$, which contradicts $\cal{F}=\{\emptyset, M_5\}\cup \cal{G}_2\cup \cal{G}_3$. Hence this case can not happen.

\item[(2.1.2.3)] $n= 4$. Now either $\cal{G}_2=\{\{1,2\}, \{1,3\}, \{1,4\}, \{1,5\}\}$  or $\cal{G}_2=\{\{2,1\}, \{2,3\}, \linebreak \{2,4\}, \{2,5\}\}$. Then $\cal{F}\cap \cal{M}_4\neq \emptyset$, which contradicts $\cal{F}=\{\emptyset, M_5\}\cup \cal{G}_2\cup \cal{G}_3$. Hence this case can not happen.
\end{itemize}
\end{itemize}

\item[(2.2)] $m$ is an odd number. Now we have $G_1\cup G_2=\cdots=G_{m-2}\cup G_{m-1}=M_5$ and thus all the 5 elements in $M_5$ belong to at least half of the sets in $\{G_1,\cdots,G_{m-1}\}$. Without loss of generality, we assume that $G_m=\{1,2,3\}$.
    We have two subcases: $|\cal{G}_2|=1$ and $|\cal{G}_2|\geq 2$.

\begin{itemize}
\item[(2.2.1)] $|\cal{G}_2|=1$. Denote $\cal{G}_2=\{H\}$. Then all the elements in $H\cup G_m$ belong to at least one of the two sets $H$ and $G_m$ and thus belong to at least half of the sets in $\cal{F}$.

\item[(2.2.2)] $|\cal{G}_2|=n\geq 2$. Denote $\cal{G}_2=\{H_1,\cdots,H_n\}$. Then for any $i,j=1,\cdots,n,i\neq j$, we must have that $H_i\cup H_j\in \cal{G}_3$.  By the condition that $\cal{F}=\{\emptyset, M_5\}\cup \cal{G}_2\cup \cal{G}_3$, we know  that for any $i=1,\cdots,n$, either $H_i\subset \{1,2,3\}$ or $H_i\cup \{1,2,3\}=M_5$, i.e. $H_i=\{4,5\}$.  If for some $i=1,\cdots,n$, $H_i=\{4,5\}$, then for any $j\in \{1,\cdots,n\}\backslash\{i\}$, we have
    $$
    H_j\in\cal{J}:= \{\{4,1,\}, \{4,2\}, \{4,3\}, \{5,1\}, \{5,2\}, \{5,3\}\}.
    $$
    But for any set $J\in\cal{J}$, $J     \nsubseteq    \{1,2,3\}$ and $J\neq \{4,5\}$.
    Hence for any $i=1,\cdots,n$, we have $H_i\subset \{1,2,3\}$.     By Case (2) in Section 2.2, we know that all the 3 elements in $\{1,2,3\}$ belong to at least $|\cal{G}_2|-1$ sets in $\cal{G}_2$. Hence all the 3 elements in $\{1,2,3\}$ belong to at least half of the sets in $\cal{F}$.
\end{itemize}
\end{itemize}
\end{itemize}

\subsubsection{$T(\cal{F}\backslash (\cal{M}_2\cup \cal{M}_3))=4$ }

Now $\cal{F}=\{\emptyset, M_5\}\cup \cal{G}_2\cup \cal{G}_3\cup \cal{G}_4$, where $\cal{G}_i$ is a nonempty subset of $\cal{M}_i$ for $i=2,3,4$. We have two subcases: $|\cal{G}_4|=1$ and $|\cal{G}_4|\geq 2$.

\begin{itemize}
\item[(1)] $|\cal{G}_4|=1$. Without loss of generality, we assume that  $\cal{G}_4=\{\{1,2,3,4\}\}$.  We have two subcases: $|\cal{G}_3|=1$ and $|\cal{G}_3|\geq 2$.

\begin{itemize}
\item[(1.1)] $|\cal{G}_3|=1$. Denote $\cal{G}_3=\{G\}$.  We have two subcases: $G\subset \{1,2,3,4\}$ and $G\nsubseteq \{1,2,3,4\}$.

\begin{itemize}
\item[(1.1.1)] $G\subset \{1,2,3,4\}$. Without loss of generality, we assume that $G=\{1,2,3\}$. We have two subcases: $|\cal{G}_2|=1$ and $|\cal{G}_2|\geq 2$.

\begin{itemize}
\item[(1.1.1.1)]  $|\cal{G}_2|=1$. Denote $\cal{G}_2=\{H\}$.  If $H\subset \{1,2,3,4\}$, then all the 2 elements in $H$ belong to at least two sets among the three sets in $\cal{G}_2\cup \cal{G}_3\cup \cal{G}_4$ and thus belong to at least half of the sets in $\cal{F}$.
If $H\nsubseteq \{1,2,3,4\}$, then $H\cup \{1,2,3,4\}=M_5$ and thus all the 5 elements in $M_5$ belong to at least one of the two sets $H$ and $\{1,2,3,4\}$.  Hence all the 3 elements in $G$ belong to at least half of the sets in $\cal{F}$.

\item[(1.1.1.2)]  $|\cal{G}_2|=n\geq 2$. Denote $\cal{G}_2=\{H_1,\cdots,H_n\}$. For any $i,j=1,\cdots,n,i\neq j$, either $H_i\cup H_j=\{1,2,3,4\}$ or $H_i\cup H_j=G$.  We have the following two subcases:
\smallskip

\quad (a) There exists a 2-element subset $\{i,j\}$ of $\{1,2,\cdots,n\}$ such that $H_i\cup H_j=\{1,2,3,4\}$, which implies that $H_i\cap H_j=\emptyset$. Then by Lemma \ref{lem-1.2}, we get that there are at least two elements in $M_5$ which belong to at least half of the sets in $\cal{F}$.

\quad (b) For any $i,j=1,2,\cdots,n,i\neq j$, we have $H_i\cup H_j=G=\{1,2,3\}$. By Case (2) in Section 2.2, we know that all the 3 elements in $\{1,2,3\}$ belong to at least $|\cal{G}_2|-1$ sets in $\cal{G}_2$. Hence all the 3 elements in $\{1,2,3\}$ belong to at least half of the sets in $\cal{F}$.
\end{itemize}

\item[(1.1.2)] $G\nsubseteq \{1,2,3,4\}$. Then $G\cup \{1,2,3,4\}=M_5$ and thus all the 5 elements in $M_5$ belong to at least one of the two sets $G$ and $\{1,2,3,4\}$.  Hence it is enough to show that there exist 2 elements in $M_5$ which belong to at least half of the sets in $\cal{G}_2$. We have two subcases: $|\cal{G}_2|=1$ and $|\cal{G}_2|\geq 2$. Without loss of generality, we assume that
    $G=\{1,2,5\}$.

\begin{itemize}
\item[(1.1.2.1)] $|\cal{G}_2|=1$. Denote $\cal{G}_2=\{H\}$. Then all the 2 elements in $H$ satisfy the condition.

\item[(1.1.2.2)] $|\cal{G}_2|=n\geq 2$. Denote $\cal{G}_2=\{H_1,\cdots,H_n\}$. For any $i,j=1,\cdots,n,i\neq j$, either $H_i\cup H_j=\{1,2,3,4\}$ or $H_i\cup H_j=G=\{1,2,5\}$.  We have the following two subcases:
\smallskip

\quad (a) There exists a 2-element subset $\{i,j\}$ of $\{1,2,\cdots,n\}$ such that $H_i\cup H_j=\{1,2,3,4\}$, which implies that $H_i\cap H_j=\emptyset$. Then by Lemma \ref{lem-1.2}, we get that there are at least two elements in $M_5$ which belong to at least half of the sets in $\cal{F}$.

\quad (b) For any $i,j=1,2,\cdots,n,i\neq j$, we have $H_i\cup H_j=G=\{1,2,5\}$. By Case (2) in Section 2.2, we know that all the 3 elements in $\{1,2,5\}$ belong to at least $|\cal{G}_2|-1$ sets in $\cal{G}_2$. Hence all the 3 elements in $\{1,2,5\}$ belong to at least half of the sets in $\cal{F}$.
\end{itemize}
\end{itemize}

\item[(1.2)] $|\cal{G}_3|=m\geq 2$. Denote $\cal{G}_3=\{G_1,\cdots,G_m\}$.  For any $i,j=1,\cdots,m,i\neq j$, either $G_i\cup G_j=M_5$ or $G_i\cup G_j=\{1,2,3,4\}$.  If $G_i\cup G_j=M_5$, then all the 5 elements in $M_5$ belong to at least one of the two sets $G_i$ and $G_j$. We have the following three subcases:

\begin{itemize}

\item[(1.2.1)] $m$ is an even number and there is a permutation $(i_1,\cdots,i_m)$ of $(1,\cdots,m)$ such that
$
G_{i_1}\cup G_{i_2}=\cdots=G_{i_{m-1}}\cup G_{i_m}=M_5.
$
Then all the 5 elements in $M_5$ belong to at least half of the sets in $\cal{G}_3$. We have two subcases: $|\cal{G}_2|=1$ and $|\cal{G}_2|\geq 2$.

\begin{itemize}
\item[(1.2.1.1)] $|\cal{G}_2|=1$. Denote $\cal{G}_2=\{H\}$. Then all the elements in $H\cup \{1,2,3,4\}$ belong to at least one of the two sets $H$ and $\{1,2,3,4\}$. Hence all the elements in $H\cup \{1,2,3,4\}$ belong to at least half of the sets in $\cal{F}$.

\item[(1.2.1.2)] $|\cal{G}_2|\geq 2$. Denote $\cal{G}_2=\{H_1,\cdots,H_n\}$. For any $i,j=1,\cdots,n,i\neq j$, either $H_i\cup H_j=\{1,2,3,4\}$ or $|H_i\cup H_j|=3$.  We have the following two subcases:
\smallskip

\quad (a) There exists a 2-element subset $\{i,j\}$ of $\{1,2,\cdots,n\}$ such that $H_i\cup H_j=\{1,2,3,4\}$, which implies that $H_i\cap H_j=\emptyset$. Then by Lemma \ref{lem-1.2}, we get that there are at least two elements in $M_5$ which belong to at least half of the sets in $\cal{F}$.

\quad (b) For any $i,j=1,2,\cdots,n,i\neq j$, we have $|H_i\cup H_j|=3$.

\quad\quad (b.1) $H_1\subset \{1,2,3,4\}$.  Without loss of generality, we assume that $H_1=\{1,2\}$.  Define
\begin{eqnarray*}
&&\cal{H}_3=\{A\in \{H_2,\cdots, H_n\}: A\cup \{1,2\}=\{1,2,3\}\}\subset \{\{1,3\},\{2,3\}\},\\
&&\cal{H}_4=\{A\in \{H_2,\cdots, H_n\}: A\cup \{1,2\}=\{1,2,4\}\}\subset \{\{1,4\},\{2,4\}\},\\
&&\cal{H}_5=\{A\in \{H_2,\cdots, H_n\}: A\cup \{1,2\}=\{1,2,5\}\}\subset\{\{1,5\},\{2,5\}\}.
\end{eqnarray*}
We have the following 7 possible subcases:

\quad\quad\quad (b.1.1) $\cal{H}_3\neq \emptyset, \cal{H}_4=\cal{H}_5=\emptyset$,

\quad\quad\quad (b.1.2) $\cal{H}_4\neq \emptyset, \cal{H}_3=\cal{H}_5=\emptyset$,

\quad\quad\quad (b.1.3) $\cal{H}_5\neq \emptyset, \cal{H}_3=\cal{H}_4=\emptyset$,

\quad\quad\quad (b.1.4) $\cal{H}_3\neq \emptyset, \cal{H}_4\neq \emptyset, \cal{H}_5=\emptyset$,

\quad\quad\quad (b.1.5) $\cal{H}_3\neq \emptyset, \cal{H}_5\neq \emptyset, \cal{H}_4=\emptyset$,

\quad\quad\quad (b.1.6) $\cal{H}_4\neq \emptyset, \cal{H}_5\neq \emptyset, \cal{H}_3=\emptyset$,

\quad\quad\quad (b.1.7) $\cal{H}_3\neq \emptyset, \cal{H}_4\neq \emptyset, \cal{H}_5\neq \emptyset$.
\smallskip

\quad As to (b.1.1), by Case (2) in Section 2.2, we know that all the 3 elements in $\{1,2,3\}$ belong to at least $|\cal{G}_2|-1$ sets in $\cal{G}_2$ and thus belong to at least half of the sets in $\cal{G}_2$. In virtue of $\cal{G}_4=\{\{1,2,3,4\}\}$, we know that all the  3 elements in $\{1,2,3\}$ belong to at least half of the sets in $\cal{F}$.
\smallskip

\quad As to (b.1.2), by following the proof of (b.1.1), we get that all the 3 elements in $\{1,2,4\}$ belong to at least half of the sets in $\cal{F}$.
\smallskip

\quad As to  (b.1.3), by following the proof of (b.1.1), we get that all the 2 elements in $\{1,2\}$ belong to at least half of the sets in $\cal{F}$.
\smallskip

\quad As to (b.1.4), by the condition (ii), we get that $|\cal{H}_3|=|\cal{H}_4|=1$. Take $\cal{H}_3=\{\{1,3\}\},\cal{H}_4=\{\{1,4\}\}$ for example. Now we get that all the 4 elements in $\{1,2,3,4\}$ belong to at least 2 sets among 4 sets in  $\cal{G}_2\cup \cal{G}_4$ and thus belong to at least half of the sets in $\cal{F}$.
\smallskip

\quad As to (b.1.5), (b.1.6) and  (b.1.7), we claim that they can not happen. Otherwise, $\{1,2,3,5\}\in \cal{F}$ or $\{1,2,4,5\}\in \cal{F}$, which contradicts the fact that $\cal{G}_4=\{\{1,2,3,4\}\}$.
\bigskip

\quad\quad (b.2) $H_1\nsubseteq\{1,2,3,4\}$. Then $H_1\cup\{1,2,3,4\}=M_5$. Without loss of generality, we assume that $H_1=\{1,5\}$.  Define
\begin{eqnarray*}
&&\cal{J}_2=\{A\in \{H_2,\cdots, H_n\}: A\cup \{1,5\}=\{1,2,5\}\}\subset \{\{1,2\},\{2,5\}\},\\
&&\cal{J}_3=\{A\in \{H_2,\cdots, H_n\}: A\cup \{1,5\}=\{1,3,5\}\}\subset \{\{1,3\},\{3,5\}\},\\
&&\cal{J}_4=\{A\in \{H_2,\cdots, H_n\}: A\cup \{1,5\}=\{1,4,5\}\}\subset\{\{1,4\},\{4,5\}\}.
\end{eqnarray*}
We have the following 7 possible subcases:

\quad\quad\quad (b.2.1) $\cal{J}_2\neq \emptyset, \cal{J}_3=\cal{J}_4=\emptyset$,

\quad\quad\quad (b.2.2) $\cal{J}_3\neq \emptyset, \cal{J}_2=\cal{J}_4=\emptyset$,

\quad\quad\quad (b.2.3) $\cal{J}_4\neq \emptyset, \cal{J}_2=\cal{J}_3=\emptyset$,

\quad\quad\quad (b.2.4) $\cal{J}_2\neq \emptyset, \cal{J}_3\neq \emptyset, \cal{J}_4=\emptyset$,

\quad\quad\quad (b.2.5) $\cal{J}_2\neq \emptyset, \cal{J}_4\neq \emptyset, \cal{J}_3=\emptyset$,

\quad\quad\quad (b.2.6) $\cal{J}_3\neq \emptyset, \cal{J}_4\neq \emptyset, \cal{J}_2=\emptyset$,

\quad\quad\quad (b.2.7) $\cal{J}_2\neq \emptyset, \cal{J}_3\neq \emptyset, \cal{J}_4\neq \emptyset$.
\smallskip

\quad As to (b.2.1),  by Case (2) in Section 2.2, we know that all the 3 elements in $\{1,2,5\}$ belong to at least $|\cal{G}_2|-1$ sets in $\cal{G}_2$ and thus belong to at least half of the sets in $\cal{G}_2$. In virtue of $\cal{G}_4=\{\{1,2,3,4\}\}$, we know that all the 2 elements in $\{1,2\}$ belong to at least half of the sets in $\cal{F}$.
\smallskip

\quad As to (b.2.2), by following the proof of (b.2.1),  we know that all the 2 elements in $\{1,3\}$ belong to at least half of the sets in $\cal{F}$.
\smallskip

\quad  As to (b.2.3), by following the proof of (b.2.1),  we know that all the 2 elements in $\{1,4\}$ belong to at least half of the sets in $\cal{F}$.
\smallskip

\quad As to (b.2.4)-(b.2.7), we claim that they can not happen. For example, as to (b.2.4),   in virtue of $\cal{G}_4=\{\{1,2,3,4\}\}$, we know that $\cal{J}_2=\{\{1,2\}\}, \cal{J}_3=\{\{1,3\}\}$. But now $\{1,5\}\cup \{1,2\}\cup \{1,3\}=\{1,2,3,5\}$, which contradicts $\cal{G}_4=\{\{1,2,3,4\}\}$.
\end{itemize}

\item[(1.2.2)] $m$ is an odd number and there is a permutation $(i_1,\cdots,i_m)$ of $(1,\cdots,m)$ such that
$
G_{i_1}\cup G_{i_2}=\cdots=G_{i_{m-2}}\cup G_{i_{m-1}}=M_5.
$
Then all the 5 elements in $M_5$ belong to at least half of the sets in $\{G_{i_1},\cdots,G_{i_{m-1}}\}$. We have two subcases: $G_{i_m}\subset \{1,2,3,4\}$ and $G_{i_m}\nsubseteq \{1,2,3,4\}$.

\begin{itemize}
\item[(1.2.2.1)] $G_{i_m}\subset \{1,2,3,4\}$. We have two subcases: $|\cal{G}_2|=1$ and $|\cal{G}_2|\geq 2$.
\smallskip

\quad (a) $|\cal{G}_2|=1$. Denote $\cal{G}_2=\{H\}$. If $H\subset \{1,2,3,4\}$, then the 2 elements in $H$ belong to at least 2 sets among three sets in $\cal{G}_2\cup \{G_{i_m}\}\cup \cal{G}_4$ and thus belong to at least  half of the sets in $\cal{F}$. If $H\nsubseteq \{1,2,3,4\}$, then $H\cup \{1,2,3,4\}=M_5$ and thus all the 3 elements in $G_{i_m}$ belong to at least half of the sets in $\cal{F}$.
\smallskip

\quad (b) $|\cal{G}_2|=n\geq 2$. Denote $\cal{G}_2=\{H_1,\cdots,H_n\}$. For any $i,j=1,\cdots,n,i\neq j$, either $H_i\cup H_j=\{1,2,3,4\}$ or $|H_i\cup H_j|=3$.  We have the following two subcases:

\smallskip

\quad\quad (b.1) There exists a 2-element subset $\{i,j\}$ of $\{1,2,\cdots,n\}$ such that $H_i\cup H_j=\{1,2,3,4\}$, which implies that $H_i\cap H_j=\emptyset$. Then by Lemma \ref{lem-1.2}, we get that there are at least two elements in $M_5$ which belong to at least half of the sets in $\cal{F}$.

\quad\quad (b.2) For any $i,j=1,2,\cdots,n,i\neq j$, we have $|H_i\cup H_j|=3$.

\quad\quad\quad  (b.2.1) $H_1\subset \{1,2,3,4\}$.  Without loss of generality, we assume that $H_1=\{1,2\}$.  Define
\begin{eqnarray*}
&&\cal{H}_3=\{A\in \{H_2,\cdots, H_n\}: A\cup \{1,2\}=\{1,2,3\}\}\subset \{\{1,3\},\{2,3\}\},\\
&&\cal{H}_4=\{A\in \{H_2,\cdots, H_n\}: A\cup \{1,2\}=\{1,2,4\}\}\subset \{\{1,4\},\{2,4\}\},\\
&&\cal{H}_5=\{A\in \{H_2,\cdots, H_n\}: A\cup \{1,2\}=\{1,2,5\}\}\subset\{\{1,5\},\{2,5\}\}.
\end{eqnarray*}
We have the following 7 possible subcases:

\quad\quad\quad (b.2.1.1) $\cal{H}_3\neq \emptyset, \cal{H}_4=\cal{H}_5=\emptyset$,

\quad\quad\quad (b.2.1.2) $\cal{H}_4\neq \emptyset, \cal{H}_3=\cal{H}_5=\emptyset$,

\quad\quad\quad (b.2.1.3) $\cal{H}_5\neq \emptyset, \cal{H}_3=\cal{H}_4=\emptyset$,

\quad\quad\quad (b.2.1.4) $\cal{H}_3\neq \emptyset, \cal{H}_4\neq \emptyset, \cal{H}_5=\emptyset$,

\quad\quad\quad (b.2.1.5) $\cal{H}_3\neq \emptyset, \cal{H}_5\neq \emptyset, \cal{H}_4=\emptyset$,

\quad\quad\quad (b.2.1.6) $\cal{H}_4\neq \emptyset, \cal{H}_5\neq \emptyset, \cal{H}_3=\emptyset$,

\quad\quad\quad (b.2.1.7) $\cal{H}_3\neq \emptyset, \cal{H}_4\neq \emptyset, \cal{H}_5\neq \emptyset$.
\smallskip

\quad As to (b.2.1.1), by Case (2) in Section 2.2, we know that all the 3 elements in $\{1,2,3\}$ belong to at least $|\cal{G}_2|-1$ sets in $\cal{G}_2$ and thus belong to at least half of the sets in $\cal{G}_2$. Note that all the 4 elements in $\{1,2,3,4\}$ belong to at least one of the two sets $G_{i_m}$ and $\{1,2,3,4\}$. Then we get that all the  3 elements in $\{1,2,3\}$ belong to at least half of the sets in $\cal{F}$.
\smallskip

\quad As to (b.2.1.2), by following the proof in (b.2.1.1), we get that all the  3 elements in $\{1,2,4\}$ belong to at least half of the sets in $\cal{F}$.
\smallskip

\quad As to  (b.2.1.3), by following the proof of (b.2.1.1), we get that the 2 elements in $\{1,2\}$ belong to at least half of the sets in $\cal{F}$.
\smallskip

\quad As to (b.2.1.4), by the condition (ii), we get that $|\cal{H}_3|=|\cal{H}_4|=1$. Take $\cal{H}_3=\{\{1,3\}\},\cal{H}_4=\{\{1,4\}\}$ for example. Now we get that all the 3 elements in $G_{i_m}$ belong to at least three sets among five sets in  $\cal{G}_2\cup \{G_{i_m}\}\cup \cal{G}_4$ and thus belong to at least half of the sets in $\cal{F}$.
\smallskip

\quad As to (b.2.1.5), (b.2.1.6) and  (b.2.1.7), we claim that they can not happen. Otherwise, $\{1,2,3,5\}\in \cal{F}$ or $\{1,2,4,5\}\in \cal{F}$, which contradicts the fact that $\cal{G}_4=\{\{1,2,3,4\}\}$.

\item[(1.2.2.2)] $G_{i_m}\nsubseteq \{1,2,3,4\}$. Then $G_{i_m}\cup \{1,2,3,4\}=M_5$ and thus it is enough to show that there exist 2 elements in $M_5$ which belong to at least half of the sets in $\cal{G}_2$.  We have two subcases: $|\cal{G}_2|=1$ and $|\cal{G}_2|\geq 2$.
\smallskip

\quad (a) $|\cal{G}_2|=1$.  Denote $\cal{G}_2=\{H\}$. Then the 2 elements in $H$ satisfy the condition.
\smallskip

\quad (b)  $|\cal{G}_2|\geq 2$. Denote $\cal{G}_2=\{H_1,\cdots,H_n\}$. For any $i,j=1,\cdots,n,i\neq j$, either $H_i\cup H_j=\{1,2,3,4\}$ or $|H_i\cup H_j|=3$.  We have the following two subcases:

\smallskip

\quad\quad (b.1) There exists a 2-element subset $\{i,j\}$ of $\{1,2,\cdots,n\}$ such that $H_i\cup H_j=\{1,2,3,4\}$, which implies that $H_i\cap H_j=\emptyset$. Then by Lemma \ref{lem-1.2}, we get that there are at least two elements in $M_5$ which belong to at least half of the sets in $\cal{F}$.

\quad\quad (b.2) For any $i,j=1,2,\cdots,n,i\neq j$, we have $|H_i\cup H_j|=3$.

\quad\quad\quad  (b.2.1) $H_1\subset \{1,2,3,4\}$.  Without loss of generality, we assume that $H_1=\{1,2\}$.  Define
\begin{eqnarray*}
&&\cal{H}_3=\{A\in \{H_2,\cdots, H_n\}: A\cup \{1,2\}=\{1,2,3\}\}\subset \{\{1,3\},\{2,3\}\},\\
&&\cal{H}_4=\{A\in \{H_2,\cdots, H_n\}: A\cup \{1,2\}=\{1,2,4\}\}\subset \{\{1,4\},\{2,4\}\},\\
&&\cal{H}_5=\{A\in \{H_2,\cdots, H_n\}: A\cup \{1,2\}=\{1,2,5\}\}\subset\{\{1,5\},\{2,5\}\}.
\end{eqnarray*}
We have the following 7 possible subcases:

\quad\quad\quad (b.2.1.1) $\cal{H}_3\neq \emptyset, \cal{H}_4=\cal{H}_5=\emptyset$,

\quad\quad\quad (b.2.1.2) $\cal{H}_4\neq \emptyset, \cal{H}_3=\cal{H}_5=\emptyset$,

\quad\quad\quad (b.2.1.3) $\cal{H}_5\neq \emptyset, \cal{H}_3=\cal{H}_4=\emptyset$,

\quad\quad\quad (b.2.1.4) $\cal{H}_3\neq \emptyset, \cal{H}_4\neq \emptyset, \cal{H}_5=\emptyset$,

\quad\quad\quad (b.2.1.5) $\cal{H}_3\neq \emptyset, \cal{H}_5\neq \emptyset, \cal{H}_4=\emptyset$,

\quad\quad\quad (b.2.1.6) $\cal{H}_4\neq \emptyset, \cal{H}_5\neq \emptyset, \cal{H}_3=\emptyset$,

\quad\quad\quad (b.2.1.7)$\cal{H}_3\neq \emptyset, \cal{H}_4\neq \emptyset, \cal{H}_5\neq \emptyset$.
\smallskip

\quad As to (b.2.1.1), by Case (2) in Section 2.2, we know that all the 3 elements in $\{1,2,3\}$ belong to at least $|\cal{G}_2|-1$ sets in $\cal{G}_2$ and thus belong to at least half of the sets in $\cal{G}_2$.
\smallskip

\quad As to (b.2.1.2), by following the proof of (b.2.1.1), we get that all the 3 elements in $\{1,2,4\}$ belong to at least half of the sets in $\cal{G}_2$.
\smallskip

\quad As to (b.2.1.3), by following the proof of (b.2.1.1), we get that all the 3 elements in $\{1,2,5\}$ belong to at least half of the sets in $\cal{G}_2$.
\smallskip

\quad As to (b.2.1.4), by the condition (ii), we get that $|\cal{H}_3|=|\cal{H}_4|=1$. Take $\cal{H}_3=\{\{1,3\}\},\cal{H}_4=\{\{1,4\}\}$ for example. Now $\cal{G}_2=\{\{1,2\},\{1,3\},\{1,4\}\}$.
By $G_{i_m}\cup \{1,2,3,4\}=M_5$, we know that $G_{i_m}=\{i,j,5\}$, where $\{i,j\}$ is a 2-element subset of $\{1,2,3,4\}$.
It follows that  the 2 elements in $\{i,j\}$ belong to at least 3 sets among  5 sets in $\cal{G}_2\cup \{G_{i_m}\}\cup \cal{G}_4 $. Hence in this case, the 2 elements in $\{i,j\}$ belong to at least half of the sets in $\cal{F}$.
\smallskip

\quad As to (b.2.1.5)-(b.2.1.7), we claim that they can not happen. Otherwise, $\{1,2,3,5\}\in \cal{F}$ or $\{1,2,4,5\}\in \cal{F}$, which contradicts the fact that $\cal{G}_4=\{\{1,2,3,4\}\}$.

\bigskip

\quad\quad\quad  (b.2.2) $H_1\nsubseteq \{1,2,3,4\}$. Without loss of generality, we assume that $H_1=\{1,5\}$. By following the proof in (1.2.1.2)(b.2), we obtain that there exist 2 elements which belong to at least half of the sets in $\cal{F}$. We omit the details.

\end{itemize}

\item[(1.2.3)] We can decompose $\{G_1,\cdots, G_m\}$ into two disjoint parts $\{G_{i_1}, \cdots, G_{i_{2k}}\}$   and\linebreak
$\{G_{i_{2k+1}},\cdots, G_{i_m}\}$, where $\{i_1,\cdots,i_m\}=\{1,\cdots,m\},m-2k\geq  2$ and

\quad\quad (i) $G_{i_1}\cup G_{i_2}=\cdots=G_{i_{2k-1}}\cup G_{i_{2k}}=M_5$,

\quad\quad (ii) for any two different indexes $i, j$ from $\{i_{2k+1}, \cdots, i_n\}, G_i\cup G_j=\{1,2,3,4\}$. \\
Then all the 5 elements in $M_5$ belong to at least half of the sets in $\{G_{i_1}, \cdots, G_{i_{2k}}\}$. By Case (2) in Section 2.2, we know that all the 4 elements in $\{1,2,3,4\}$ belong to at least $|\{G_{i_{2k+1}},\cdots, G_{i_m}\}|-1$ sets in $\{G_{i_{2k+1}},\cdots, G_{i_m}\}$ and thus belong to at least half of the sets in $\{G_{i_{2k+1}},\cdots, G_{i_m}\}$.  Hence all the 4 elements in $\{1,2,3,4\}$ belong to at least half of the sets in $\cal{G}_3$ and thus  it is enough to show that there exist 2 elements in $\{1,2,3,4\}$ which belong to at least half of the sets in $\cal{G}_2\cup \cal{G}_4$. We have two subcases: $|\cal{G}_2|=1$ and $|\cal{G}_2|\geq 2$.
\smallskip

\quad (a) $|\cal{G}_2|=1$.  Denote $\cal{G}_2=\{H\}$. Then all the 4 elements in $\{1,2,3,4\}$ belong to at least half of the sets in $\cal{G}_2\cup \cal{G}_4$.
\smallskip

\quad (b)  $|\cal{G}_2|\geq 2$. Denote $\cal{G}_2=\{H_1,\cdots,H_n\}$. For any $i,j=1,\cdots,n,i\neq j$, either $H_i\cup H_j=\{1,2,3,4\}$ or $|H_i\cup H_j|=3$.  We have the following two subcases:

\smallskip

\quad\quad (b.1) There is  a 2-element subset $\{i,j\}$ of $\{1,2,\cdots,n\}$ such that $H_i\cup H_j=\{1,2,3,4\}$, which implies that $H_i\cap H_j=\emptyset$. Then by Lemma \ref{lem-1.2}, we get that there are at least two elements in $M_5$ which belong to at least half of the sets in $\cal{F}$.

\quad\quad (b.2) For any $i,j=1,2,\cdots,n,i\neq j$, we have $|H_i\cup H_j|=3$.

Now by following the proof  in (1.2.1.2)(b), we can get that there exist 2 elements in $\{1,2,3,4\}$ which belong to at least half of the sets in $\cal{G}_2\cup \cal{G}_4$.
\end{itemize}
\end{itemize}

\item[(2)] $|\cal{G}_4|\geq 2$. By Case (2) in Section 2.2, we know that all the 5 elements in $M_5$ belong to at least $|\cal{G}_4|-1$ sets in $\cal{G}_4$ and thus belong to at least half of the sets in $\cal{G}_4$. Hence it is enough to show that there exist 2 elements in $M_5$ which belong to at least half of the sets in $\cal{G}_2\cup \cal{G}_3$. We have two subcases: $|\cal{G}_3|=1$ and $|\cal{G}_3|\geq 2$.

\begin{itemize}
\item[(2.1)] $|\cal{G}_3|=1$. Without loss of generality, we assume that  $\cal{G}_3=\{\{1,2,3\}\}$. We have two subcases: $|\cal{G}_2|=1$ and $|\cal{G}_2|\geq 2$.

\begin{itemize}

\item[(2.1.1)] $|\cal{G}_2|=1$. Denote $\cal{G}_2=\{H\}$. Then all the elements in $H\cup \{1,2,3\}$ belong to at least one of the two sets  $H$ and $\{1,2,3\}$ and thus belong to at least half of the sets in $\cal{G}_2\cup \cal{G}_3$.

\item[(2.1.2)] $|\cal{G}_2|=n\geq 2$. Denote $\cal{G}_2=\{H_1,\cdots,H_n\}$. For any $i,j=1,\cdots,n,i\neq j$, either $H_i\cup H_j\in \cal{G}_4$ or $H_i\cup H_j=\{1,2,3\}$. We have two subcases:

\begin{itemize}
\item[(2.1.2.1)]  There exists a 2-element subset of $\{1,\cdots,n\}$ such that $H_i\cup H_j\in \cal{G}_4$, which implies that $H_i\cap H_j=\emptyset$. Then by Lemma \ref{lem-1.2}, we get that there are at least two elements in $M_5$ which belong to at least half of the sets in $\cal{F}$.

\item[(2.1.2.2)]  For any $i,j=2,\cdots,n,i\neq j$, $H_i\cup H_j=\{1,2,3\}$.   Now by Case (2) in Section 2.2, we know that all the 3 elements in $\{1,2,3\}$ belong to at least $|\cal{G}_2|-1$ sets in $\cal{G}_2$ and thus all the 3 elements in $\{1,2,3\}$ belong to at least half of the  sets in $\cal{G}_2\cup \cal{G}_3$.
\end{itemize}

\end{itemize}

\item[(2.2)] $|\cal{G}_3|=m\geq 2$. Denote $\cal{G}_3=\{G_1,\cdots,G_m\}$.  For any $i,j=1,\cdots,m,i\neq j$, either $G_i\cup G_j=M_5$ or $|G_i\cup G_j|=4$.  If $G_i\cup G_j=M_5$, then all the 5 elements in $M_5$ belong to at least one of the two sets $G_i$ and $G_j$. We have the following three subcases:

\begin{itemize}
\item[(2.2.1)]  $m$ is an even number and there is a permutation $(i_1,\cdots,i_m)$ of $(1,\cdots,m)$ such that
$
G_{i_1}\cup G_{i_2}=\cdots=G_{i_{m-1}}\cup G_{i_m}=M_5.
$
Then all the 5 elements in $M_5$ belong to at least half of the sets in $\cal{G}_3$. Hence it is enough to show that there exist 2 elements in $M_5$ which belong to at least half of the sets in $\cal{G}_2$ or in $\cal{G}_2\cup \cal{G}_4$.   We have two subcases: $|\cal{G}_2|=1$ and $|\cal{G}_2|\geq 2$.

\begin{itemize}
\item[(2.2.1.1)] $|\cal{G}_2|=1$. Denote $\cal{G}_2=\{H\}$. Then the 2 elements in $H$ satisfy the condition.

 \item[(2.2.1.2)] $|\cal{G}_2|\geq 2$. Denote $\cal{G}_2=\{H_1,\cdots,H_n\}$. For any $i,j=1,\cdots,n,i\neq j$, either $H_i\cup H_j\in \cal{G}_4$ or $H_i\cup H_j\in \cal{G}_3$.  
     We have the following two subcases:
     \smallskip

  \quad (a)  There exists a 2-element subset of $\{1,\cdots,n\}$ such that $H_i\cup H_j\in \cal{G}_4$, which implies that $H_i\cap H_j=\emptyset$. Then by Lemma \ref{lem-1.2}, we get that there are at least two elements in $M_5$ which belong to at least half of the sets in $\cal{F}$.
       \smallskip

  \quad (b) For any $i,j=1,\cdots,n,i\neq j$,  $H_i\cup H_j\in \cal{G}_3$.  Now we can check that $n\leq 4$.  And so we have only three subcases, which are $n=2,n=3,n=4$.
     \smallskip

\quad\quad (b.1) $n=2.$ Without loss of generality, we assume $\cal{G}_2=\{\{1,2\},\{1,3\}\}$. Now all the 3 elements in $\{1,2,3\}$ belong to at least half of the sets in $\cal{G}_2$.
 \smallskip

\quad\quad (b.2) $n=3.$  Now we have two subcases:   $\cap_{i=1}^3H_i=\emptyset$ and $|\cap_{i=1}^3H_i|=1$.
 \smallskip

\quad\quad\quad  (b.2.1) $\cap_{i=1}^3H_i=\emptyset$. Without loss of generality, we assume that $\cal{G}_2=\{\{1,2\},\{1,3\},\{2,3\}\}$.   Now all the 3 elements in $\{1,2,3\}$ belong to at least half of the sets in $\cal{G}_2$.
 \smallskip

\quad\quad\quad  (b.2.2) $|\cap_{i=1}^3H_i|=1$.  Without loss of generality, we assume that $\cal{G}_2=\{\{1,2\},\{1,3\},\{1,4\}\}$.
Notice that $\{1,2,3,4\}\in \cal{G}_4\subset \cal{M}_4=\{\{1,2,3,4\},\linebreak  \{1,2,3,5\}, \{1,2,4,5\}, \{1,3,4,5\}, \{2,3,4,5\}\}$ and $|\cal{G}_4|\geq 2$.  If $\cal{G}_4=\cal{M}_4$, then all 5 elements in $M_5$ belong to at least half of the sets in $\cal{G}_2\cup \cal{G}_4$.  If $2\leq |\cal{G}_4|\leq 4$, we have the following 14 cases:

\quad\quad\quad\quad   (b.2.2.1) $\cal{G}_4=\{\{1,2,3,4\}, \{1,2,3,5\}\}$.

\quad\quad\quad\quad   (b.2.2.2) $\cal{G}_4=\{\{1,2,3,4\}, \{1,2,4,5\}\}$.

\quad\quad\quad\quad   (b.2.2.3) $\cal{G}_4=\{\{1,2,3,4\}, \{1,3,4,5\}\}$.

\quad\quad\quad\quad   (b.2.2.4) $\cal{G}_4=\{\{1,2,3,4\}, \{2,3,4,5\}\}$.

\quad\quad\quad\quad   (b.2.2.5) $\cal{G}_4=\{\{1,2,3,4\}, \{1,2,3,5\},\{1,2,4,5\}\}$.

\quad\quad\quad\quad   (b.2.2.6) $\cal{G}_4=\{\{1,2,3,4\}, \{1,2,3,5\},\{1,3,4,5\}\}$.

\quad\quad\quad\quad   (b.2.2.7) $\cal{G}_4=\{\{1,2,3,4\}, \{1,2,3,5\},\{2,3,4,5\}\}$.

\quad\quad\quad\quad   (b.2.2.8) $\cal{G}_4=\{\{1,2,3,4\}, \{1,2,4,5\},\{1,3,4,5\}\}$.

\quad\quad\quad\quad   (b.2.2.9) $\cal{G}_4=\{\{1,2,3,4\}, \{1,2,4,5\},\{2,3,4,5\}\}$.

\quad\quad\quad\quad   (b.2.2.10) $\cal{G}_4=\{\{1,2,3,4\}, \{1,3,4,5\},\{2,3,4,5\}\}$.

\quad\quad\quad\quad   (b.2.2.11) $\cal{G}_4=\{\{1,2,3,4\}, \{1,2,3,5\},\{1,2,4,5\},\{1,3,4,5\}\}$.

\quad\quad\quad\quad   (b.2.2.12) $\cal{G}_4=\{\{1,2,3,4\}, \{1,2,3,5\},\{1,2,4,5\},\{2,3,4,5\}\}$.

\quad\quad\quad\quad   (b.2.2.13) $\cal{G}_4=\{\{1,2,3,4\}, \{1,2,3,5\},\{1,3,4,5\},\{2,3,4,5\}\}$.

\quad\quad\quad\quad   (b.2.2.14) $\cal{G}_4=\{\{1,2,3,4\}, \{1,2,4,5\},\{1,3,4,5\},\{2,3,4,5\}\}$.
\smallskip

 For the above 14 cases, we can easily check that there exist at least 2 elements in $M_5$ which belong to at least half of the sets in $\cal{G}_2\cup \cal{G}_4$.

 \smallskip

\quad\quad (b.3) $n=4.$  Now, without loss of generality, we assume that $\cal{G}_2=\{\{1,2\},\linebreak \{1,3\},\{1,4\}, \{1,5\}\}$. Now $\{\{1,2,3,4\}, \{1,2,3,5\}, \{1,2,4,5\},\{1,3,4,5\}\}\subset \cal{G}_4$.  So we have the following two subcases:

\quad\quad\quad (b.3.1) $\cal{G}_4=\cal{M}_4$.

\quad\quad\quad (b.3.2) $\cal{G}_4=\{\{1,2,3,4\}, \{1,2,3,5\}, \{1,2,4,5\},\{1,3,4,5\}\}$.

For these two cases, we can check  that that there exists at least 2 elements in $M_5$ which belong to at least half of the sets in $\cal{G}_2\cup \cal{G}_4$.
\end{itemize}

\item[(2.2.2)]  $m$ is an odd number and there is a permutation $(i_1,\cdots,i_m)$ of $(1,\cdots,m)$ such that
\begin{eqnarray}\label{2.2.2-a}
G_{i_1}\cup G_{i_2}=\cdots=G_{i_{m-2}}\cup G_{i_{m-1}}=M_5.
\end{eqnarray}
Then all the 5 elements in $M_5$ belong to at least half of the sets in $\{G_{i_1},\cdots,G_{i_{m-1}}\}$. Hence it is enough to show that there exist 2 elements in $M_5$ which belong to at least half of the sets in $\cal{G}_2\cup \{G_{i_m}\}$ or $\cal{G}_2\cup \{G_{i_m}\}\cup \cal{G}_4$.  Without loss of generality, we assume that $G_{i_m}=\{1,2,3\}$.   We have two subcases: $|\cal{G}_2|=1$ and $|\cal{G}_2|\geq 2$.

\begin{itemize}
\item[(2.2.2.1)] $|\cal{G}_2|=1$.  Denote $\cal{G}_2=\{H\}$. The all the elements in $H\cup G_{i_m}$ belong to at least one of the two sets in $\cal{G}_2\cup \{G_{i_m}\}$.

\item[(2.2.2.2)] $|\cal{G}_2|=n\geq 2$. Denote $\cal{G}_2=\{H_1,\cdots,H_n\}$. For any $i,j=1,\cdots,n,i\neq j$, either $H_i\cup H_j\in \cal{G}_4$ or $H_i\cup H_j\in \cal{G}_3$.  
     We have the following two subcases:
     \smallskip

  \quad (a)  There exists a 2-element subset of $\{1,\cdots,n\}$ such that $H_i\cup H_j\in \cal{G}_4$, which implies that $H_i\cap H_j=\emptyset$. Then by Lemma \ref{lem-1.2}, we get that there are at least two elements in $M_5$ which belong to at least half of the sets in $\cal{F}$.
       \smallskip

        \quad (b) For any $i,j=1,\cdots,n,i\neq j$,  $H_i\cup H_j\in \cal{G}_3$.  Now we can check that $n\leq 4$.  And so we have only three subcases, which are $n=2,n=3,n=4$.
     \smallskip

\quad\quad (b.1) $n=2$. Now $\cal{G}_2=\{H_1,H_2\}$. We have two subcases:
  \smallskip

\quad\quad\quad (b.1.1) There exists $j\in \{1,2\}$ such that $H_j\subset G_{i_m}$.  Then all the 2 elements in $H_j$ belong to at least 2 sets among the three sets in $\cal{G}_2\cup \{G_{i_m}\}$.

 \smallskip

\quad\quad\quad (b.1.2) For any $j=1,2$, $H_j\nsubseteqq G_{i_m}$.  We have the following two subcases:

\smallskip

\quad\quad\quad\quad (b.1.2.1) There exists $j\in \{1,2\}$ such that $H_j=\{4,5\}$. Without loss of generality, we assume that $H_1=\{4,5\}$.  Then $H_2\in \{\{1,4\}, \{2,4\}, \{3,4\}, \linebreak \{1,5\}, \{2,5\}, \{3,5\}\}$. Take $H_2=\{1,4\}$ for example. Then
all the 2 elements in $\{1,4\}$ belong to at least 2 sets among the three sets in $\cal{G}_2\cup \{G_{i_m}\}$.
\smallskip

\quad\quad\quad\quad (b.1.2.2)  For any $j=1,2$, $H_j\neq \{4,5\}$. Without loss of generality, we assume that $H_1=\{1,4\}$. Then
$H_2\in \{\{2,4\}, \{3,4\}, \{1,5\}\}$.

\quad If $H_2\in \{\{2,4\}, \{3,4\}\}$, one can easily check that there at least 2 elements in $M_5$ which belong to at least half of the sets in $\cal{G}_2\cup \{G_{i_m}\}$.

\quad If $H_2=\{1,5\}$. Then $H_1\cup H_2=\{1,4,5\}\in \{G_{i_1},\cdots,G_{i_{m-1}}\}$.  Without loss of generality, we assume that
$G_{i_1}=\{1,4,5\}$. By $G_{i_1}\cup G_{i_2}=M_5$ and the assumption that $G_{i_m}=\{1,2,3\}$, we get that  $G_{i_2}\in \{\{2,3,4\}, \{2,3,5\}\}$.  Now we exchange $G_{i_2}$ and $G_{i_m}$, i.e. we rewrite (\ref{2.2.2-a}) as follows:
\begin{eqnarray*}
G_{i_1}\cup G_{i_m}=G_{i_3}\cup G_{i_4}=\cdots=G_{i_{m-2}}\cup G_{i_{m-1}}=M_5,
\end{eqnarray*}
and it is enough to show that there exist 2 elements in $M_5$ which belong to at least half of the sets in $\cal{G}_2\cup \{G_{i_2}\}$, where $\cal{G}_2\cup \{G_{i_2}\}=\{\{1,4\},\{1,5\},\{2,3,4\}\}$ or $\{\{1,4\},\{1,5\},\{2,3,5\}\}$. For these two cases, we can easily check that  there exist 2 elements in $M_5$ which belong to at least half of the sets in $\cal{G}_2\cup \{G_{i_2}\}$.

 \smallskip

\quad\quad (b.2) $n=3$.  Now $\cal{G}_2=\{H_1,H_2, H_3\}$. We have two subcases:
\smallskip

\quad\quad\quad (b.2.1) There exists $j\in \{1,2,3\}$ such that $H_j\subset G_{i_m}$.  Then all the 2 elements in $H_j$ belong to at least 2 sets among the 4 sets in $\cal{G}_2\cup \{G_{i_m}\}$.

 \smallskip

\quad\quad\quad (b.2.2) For any $j=1,2,3$, $H_j\nsubseteqq G_{i_m}$.  We have the following two subcases:

\smallskip

\quad\quad\quad\quad (b.2.2.1) There exists $j\in \{1,2,3\}$ such that $H_j=\{4,5\}$. Without loss of generality, we assume that $H_1=\{4,5\}$.  Then $\{H_2,H_3\}\subset  \{\{1,4\}, \{2,4\},\linebreak \{3,4\},  \{1,5\}, \{2,5\}, \{3,5\}\}$. Take $H_2=\{1,4\}$ for example. Then all the 2 elements in $\{1,4\}$ belong to at least 2 sets among the 4 sets in $\cal{G}_2\cup \{G_{i_m}\}$.
\smallskip

\quad\quad\quad\quad (b.2.2.2)  For any $j=1,2,3$, $H_j\neq \{4,5\}$. Without loss of generality, we assume that $H_1=\{1,4\}$. Then
$\{H_2,H_3\}\subset  \{\{2,4\}, \{3,4\}, \{1,5\}\}$.  Obviously, at least one of the two sets in $\{H_2,H_3\}$ belongs to $\{\{2,4\}, \{3,4\}\}$ and then we can easily check that  there exist at least 2 elements in $M_5$ which belong to at least two sets in the 4 sets in $\cal{G}_2\cup \{G_{i_m}\}$.

 \smallskip

\quad\quad (b.3) $n=4$. Now $\cal{G}_2=\{H_1,H_2, H_3,H_4\}$ and  we have only the following 5 subcases:

\quad\quad\quad  (b.3.1) $\cal{G}_2=\{\{1,2\},\{1,3\}, \{1,4\},\{1,5\}\}$.

\quad\quad\quad  (b.3.2) $\cal{G}_2=\{\{2,1\},\{2,3\}, \{2,4\},\{2,5\}\}$.

\quad\quad\quad  (b.3.3) $\cal{G}_2=\{\{3,1\},\{3,2\}, \{3,4\},\{3,5\}\}$.

\quad\quad\quad  (b.3.4) $\cal{G}_2=\{\{4,1\},\{4,2\}, \{4,3\},\{4,5\}\}$.

\quad\quad\quad  (b.3.5) $\cal{G}_2=\{\{5,1\},\{5,2\}, \{5,3\},\{5,4\}\}$.

\smallskip

\quad As to (b.3.1), we know that $\{\{1,2,3,4\},\{1,2,3,5\}, \{1,2,4,5\}, \{1,3,4,5\}\}\subset \cal{G}_4$. If $\cal{G}_4=\{\{1,2,3,4\},\{1,2,3,5\}, \{1,2,4,5\}, \{1,3,4,5\}\}$, then all the 3 elements in $\{1,2,3\}$ belong to at least 5 sets among 9 sets in $\cal{G}_2\cup \{G_{i_m}\}\cup \cal{G}_4$ and thus belong to at least half of the sets in $\cal{F}$.  If $\cal{G}_4=\{\{1,2,3,4\},\linebreak \{1,2,3,5\}, \{1,2,4,5\}, \{1,3,4,5\},\{2,3,4,5\}\}$, then all the 3 elements in $\{1,2,3\}$ belong to at least 5 sets among 10 sets in $\cal{G}_2\cup \{G_{i_m}\}\cup \cal{G}_4$ and thus belong to at least half of the sets in $\cal{F}$.
\smallskip

\quad As to (b.3.2) - (b.3.5), by following the above analysis, we can obtain that there exist at least 2 elements in $M_5$ which belong to at least half of the sets in $\cal{F}$.

\end{itemize}

\item[(2.2.3)]  We can decompose $\{G_1,\cdots, G_m\}$ into two disjoint parts $\{G_{i_1}, \cdots, G_{i_{2k}}\}$   and\linebreak
$\{G_{i_{2k+1}},\cdots, G_{i_m}\}$, where $\{i_1,\cdots,i_m\}=\{1,\cdots,m\},m-2k\geq  2$ and

\quad\quad (i) $G_{i_1}\cup G_{i_2}=\cdots=G_{i_{2k-1}}\cup G_{i_{2k}}=M_5$,

\quad\quad (ii) for any two different indexes $i, j$ from $\{i_{2k+1}, \cdots, i_m\}, |G_i\cup G_j|=4$. \\
Then all the 5 elements in $M_5$ belong to at least half of the sets in $\{G_{i_1}, \cdots, G_{i_{2k}}\}$ and thus it is enough to show that there exist two elements in $M_5$ which belong to at least half of the sets in $\cal{G}_2\cup \{G_{i_{2k+1}},\cdots, G_{i_m}\}$ or $\cal{G}_2\cup \{G_{i_{2k+1}},\cdots, G_{i_m}\}\cup \cal{G}_4$. Without loss of generality, we assume that $G_{i_{2k+1}}=\{1,2,3\}$. Then for any $j=2k+2,\ldots, m$, we have
$$
G_{i_j}\in \{\{1,2,4\}, \{1,3,4\}, \{2,3,4\}, \{1,2,5\}, \{1,3,5\}, \{2,3,5\}\}.
$$
By $\{1,2,4\}\cup \{2,3,5\}=\{1,2,4\}\cup\{1,3,5\}=\{1,3,4\}\cup \{1,2,5\}=\{1,3,4\}\cup \{2,3,5\}=\{2,3,4\}\cup \{1,3,5\}=\{2,3,4\}\cup\{1,2,5\}=M_5,$ and the condition (ii), we get that $m-2k\leq 4$. Thus we have the following 3 subcases:

\begin{itemize}
\item[(2.2.3.1)]  $m-2k=2$. Take $\{G_{i_{2k+1}},\cdots, G_{i_m}\}=\{\{1,2,3\}, \{1,2,4\}\}$ for example. We have 2 subcases:
\smallskip

\quad (a) $|\cal{G}_2|=1$.  Denote $\cal{G}_2=\{H\}$. Now the 2 elements in $\{1,2\}$ belong to at least 2 sets among three sets $H, \{1,2,3\}$ and $\{1,2,4\}$.

\quad (b) $|\cal{G}_2|=n\geq 2$.  Denote $\cal{G}_2=\{H_1,\cdots,H_n\}$.  We have the following two subcases:
     \smallskip

  \quad\quad  (b.1)  There exists a 2-element subset of $\{1,\cdots,n\}$ such that $H_i\cup H_j\in \cal{G}_4$, which implies that $H_i\cap H_j=\emptyset$. Then by Lemma \ref{lem-1.2}, we get that there are at least two elements in $M_5$ which belong to at least half of the sets in $\cal{F}$.
       \smallskip

        \quad\quad  (b.2) For any $i,j=1,\cdots,n,i\neq j$,  $H_i\cup H_j\in \cal{G}_3$.  Now we can check that $n\leq 4$.  And so we have only three cases, which are $n=2,n=3,n=4$.
 \smallskip

\quad\quad\quad  (b.2.1) $n=2$.  Now the 2 elements in $\{1,2\}$ belong to at least 2 sets among four sets in $\cal{G}_2\cup \{G_{i_{2k+1}},\cdots, G_{i_m}\}$.
  \smallskip

  \quad\quad\quad  (b.2.2) $n=3$. Now $\cal{G}_2=\{H_1,H_2,H_3\}$.
 \smallskip

  \quad\quad\quad\quad   (b.2.2.1) There exists $j\in \{1,2,3\}$ such that $H_j=\{1,2\}$. Then all the 2 elements in $\{1,2\}$ belong to at least 3 sets among 5 sets in $\cal{G}_2\cup \{G_{i_{2k+1}},\cdots, G_{i_m}\}$.
  \smallskip

    \quad\quad\quad\quad   (b.2.2.2)  There exists $j\in \{1,2,3\}$ such that $H_j=\{4,5\}$. Without loss of generality, we assume that $H_1=\{4,5\}$ and thus $H_1\cup \{1,2,3\}=M_5$. Hence  it is enough to show that there exist two elements in $M_5$ which belong to at least half of the sets in $\{H_2,H_3\}\cup \{\{1,2,4\}\}$. Now we have
    \begin{eqnarray}\label{3.2}
    \{H_2,H_3\}\subset \{\{1,4\}, \{2,4\}, \{3,4\}, \{1,5\},\{2,5\}, \{3,5\}\}.
    \end{eqnarray}
  If $\{1,4\}\in \{H_2,H_3\}$, then the 2 elements in $\{1,4\}$ belong to at least 2 sets among three sets in $\{H_2,H_3\}\cup \{\{1,2,4\}\}$.  If $\{2,4\}\in \{H_2,H_3\}$, then the 2 elements in $\{2,4\}$ belong to at least 2 sets among three sets in $\{H_2,H_3\}\cup \{\{1,2,4\}\}$. Hence we need only consider the following 4 cases: $\{H_2,H_3\}=\{\{3,4\}, \{3,5\}\}, \{\{1,5\}, \{2,5\}\}, \{\{1,5\}, \{3,5\}\}, \{\{2,5\}, \{3,5\}\}$.  We can easily check that in all these 4 cases, there exists 2 elements in $M_5$ which belong to at least half of the sets in $\{H_2,H_3\}\cup \{\{1,2,4\}\}$.
  \smallskip

   \quad\quad\quad\quad   (b.2.2.3)  There exists $j\in \{1,2,3\}$ such that $H_j=\{3,5\}$. Similar to the analysis in (b.2.1.2), we can get that there exist two elements in $M_5$ which belong to at least half of the sets in $\cal{G}_2\cup \{G_{i_{2k+1}},\cdots, G_{i_m}\}$.

  \smallskip

   \quad\quad\quad\quad   (b.2.2.4)   $\cal{G}_2\subset \{\{1,3\}, \{1,4\}, \{1,5\}, \{2,3\}, \{2,4\},\{2,5\}, \{3,4\}\}$.  Now we have the following 6 cases:

   \quad\quad\quad\quad\quad    (b.2.2.4-1)  $\cal{G}_2=\{\{1,3\}, \{1,4\}, \{1,5\}\}$.

    \quad\quad\quad\quad\quad    (b.2.2.4-2) $\cal{G}_2=\{\{1,3\}, \{1,4\}, \{3,4\}\}$.

    \quad\quad\quad\quad\quad    (b.2.2.4-3)  $\cal{G}_2=\{\{1,3\}, \{2,3\}, \{3,4\}\}$.

   \quad\quad\quad\quad\quad    (b.2.2.4-4)  $\cal{G}_2=\{\{1,4\}, \{2,4\}, \{3,4\}\}$.

 \quad\quad\quad\quad\quad    (b.2.2.4-5)  $\cal{G}_2=\{\{2,3\}, \{2,4\}, \{3,4\}\}$.

\quad\quad\quad\quad\quad    (b.2.2.4-6)  $\cal{G}_2=\{\{2,3\}, \{2,4\}, \{2,5\}\}$.
\smallskip

\quad As to Cases (b.2.2.4-2) - (b.2.2.4-5), we  can easily check that there exist two elements in $M_5$ which belong to at least half of the sets in $\cal{G}_2\cup \{G_{i_{2k+1}},\cdots, G_{i_m}\}$.

\smallskip

\quad As to  (b.2.2.4-1), we know that $\{\{1,2,3,4\}, \{1,2,3,5\}, \{1,2,4,5\},\{1,3,4,5\}\}\linebreak \subset \cal{G}_4$.


\quad\quad If $
\cal{G}_4=\{\{1,2,3,4\}, \{1,2,3,5\}, \{1,2,4,5\}, \{1,3,4,5\}\},
$
then all the 3 elements in $\{1,3,4\}$ belong to at least half of the sets in $\cal{G}_2\cup \{G_{i_{2k+1}},\cdots, G_{i_m}\}\cup \cal{G}_4$ and thus belong to at least half of the sets in $\cal{F}$.


\quad\quad If $
\cal{G}_4=\{\{1,2,3,4\}, \{1,2,3,5\}, \{1,2,4,5\}, \{1,3,4,5\}, \{2,3,4,5\}\},
$
then all the 5 elements in $M_5$ belong to at least half of the sets in $\cal{G}_2\cup \{G_{i_{2k+1}},\cdots, G_{i_m}\}\cup \cal{G}_4$ and thus belong to at least half of the sets in $\cal{F}$.
\smallskip

\quad As to  (b.2.2.4-6), we know that $\{\{1,2,3,4\}, \{1,2,3,5\}, \{1,2,4,5\}, \{2,3,4,5\}\}\linebreak \subset \cal{G}_4$.

\quad\quad If $\cal{G}_4=\{\{1,2,3,4\}, \{1,2,3,5\}, \{1,2,4,5\}, \{2,3,4,5\}\}$, then all 4 elements in $\{1,2,3,4\}$ belong to at least half of the sets in $\cal{G}_2\cup \{G_{i_{2k+1}},\cdots, G_{i_m}\}\cup \cal{G}_4$ and thus belong to at least half of the sets in $\cal{F}$.

\quad\quad If $\cal{G}_4=\{\{1,2,3,4\}, \{1,2,3,5\}, \{1,2,4,5\}, \{2,3,4,5\},\{1,3,4,5\}\}$, then all 5 elements in $M_5$ belong to at least half of the sets in $\cal{G}_2\cup \{G_{i_{2k+1}},\cdots, G_{i_m}\}\cup \cal{G}_4$ and thus belong to at least half of the sets in $\cal{F}$.

\smallskip

\quad\quad\quad  (b.2.1) $n=4$.  Now $\cal{G}_2=\{H_1,H_2, H_3,H_4\}$ and  we have only the following 5 subcases:

\quad\quad\quad  (b.2.1.1) $\cal{G}_2=\{\{1,2\},\{1,3\}, \{1,4\},\{1,5\}\}$.

\quad\quad\quad  (b.2.1.2) $\cal{G}_2=\{\{2,1\},\{2,3\}, \{2,4\},\{2,5\}\}$.

\quad\quad\quad (b.2.1.3) $\cal{G}_2=\{\{3,1\},\{3,2\}, \{3,4\},\{3,5\}\}$.

\quad\quad\quad  (b.2.1.4) $\cal{G}_2=\{\{4,1\},\{4,2\}, \{4,3\},\{4,5\}\}$.

\quad\quad\quad (b.2.1.5) $\cal{G}_2=\{\{5,1\},\{5,2\}, \{5,3\},\{5,4\}\}$.

\smallskip

\quad As to (b.2.1.1) and (b.2.1.2), the 2 elements in $\{1,2\}$ belong to at least half of the sets in $\cal{G}_2\cup \{G_{i_{2k+1}},\cdots, G_{i_m}\}$. \smallskip

\quad As to (b.2.1.3),  the 2 elements in $\{2,3\}$ belong to at least half of the sets in $\cal{G}_2\cup \{G_{i_{2k+1}},\cdots, G_{i_m}\}$.
\smallskip

\quad As to (b.2.1.4),  the 2 elements in $\{1,4\}$ belong to at least half of the sets in $\cal{G}_2\cup \{G_{i_{2k+1}},\cdots, G_{i_m}\}$.

\smallskip

\quad As to (b.2.1.5),  the 2 elements in $\{2,5\}$ belong to at least half of the sets in $\cal{G}_2\cup \{G_{i_{2k+1}},\cdots, G_{i_m}\}$.
\bigskip

\item[(2.2.3.2)]  $m-2k=3$. Now we have the following 9 cases:

\quad (2.2.3.2-1) $\{G_{i_{2k+1}},\cdots, G_{i_m}\}=\{\{1,2,3\}, \{1,2,4\},\{ 1,3,4\} \}$.

\quad (2.2.3.2-2) $\{G_{i_{2k+1}},\cdots, G_{i_m}\}=\{\{1,2,3\}, \{1,2,4\},\{ 2,3,4\} \}$.

\quad (2.2.3.2-3) $\{G_{i_{2k+1}},\cdots, G_{i_m}\}=\{\{1,2,3\}, \{1,2,4\},\{ 1,2,5\} \}$.

\quad (2.2.3.2-4) $\{G_{i_{2k+1}},\cdots, G_{i_m}\}=\{\{1,2,3\}, \{1,3,4\},\{ 2,3,4\} \}$.

\quad (2.2.3.2-5) $\{G_{i_{2k+1}},\cdots, G_{i_m}\}=\{\{1,2,3\}, \{1,3,4\},\{ 1,3,5\} \}$.

\quad (2.2.3.2-6) $\{G_{i_{2k+1}},\cdots, G_{i_m}\}=\{\{1,2,3\}, \{2,3,4\},\{ 2,3,5\} \}$.

\quad (2.2.3.2-7) $\{G_{i_{2k+1}},\cdots, G_{i_m}\}=\{\{1,2,3\}, \{1,2,5\},\{ 1,3,5\} \}$.

\quad (2.2.3.2-8) $\{G_{i_{2k+1}},\cdots, G_{i_m}\}=\{\{1,2,3\}, \{1,2,5\},\{ 2,3,5\} \}$.

\quad (2.2.3.2-9) $\{G_{i_{2k+1}},\cdots, G_{i_m}\}=\{\{1,2,3\}, \{1,3,5\},\{ 2,3,5\} \}$.
\bigskip

\quad (2.2.3.2-1):  We have two subcases: $|\cal{G}_2|=1$ and $|\cal{G}_2|\geq 2$.

\quad\quad (a) $|\cal{G}_2|=1$. Denote $\cal{G}_2=\{H\}$. Now all the 4 elements in $\{1,2,3,4\}$ belong to at least 2 sets among the four sets in $\cal{G}_2\cup \{G_{i_{2k+1}},\cdots, G_{i_m}\}$.

\quad\quad (b) $|\cal{G}_2|=n\geq 2$.  Denote $\cal{G}_2=\{H_1,\cdots,H_n\}$.  We have the following two subcases:
     \smallskip

  \quad\quad\quad  (b.1)  There exists a 2-element subset of $\{1,\cdots,n\}$ such that $H_i\cup H_j\in \cal{G}_4$, which implies that $H_i\cap H_j=\emptyset$. Then by Lemma \ref{lem-1.2}, we get that there are at least two elements in $M_5$ which belong to at least half of the sets in $\cal{F}$.
       \smallskip

        \quad\quad\quad  (b.2) For any $i,j=1,\cdots,n,i\neq j$,  $H_i\cup H_j\in \cal{G}_3$.  Now we can check that $n\leq 4$.  And so we have only three cases, which are $n=2,n=3,n=4$.

\quad\quad\quad\quad  (b.2.1) $n=2$.  Now $\cal{G}_2=\{H_1,H_2\}$
 and thus $\cal{G}_2\cap (\cal{M}_2\backslash\{1,5\})\neq \emptyset$, where
 $$
 \cal{M}_2=\{\{1,2\}, \{1,3\}, \{1,4\}, \{1,5\}, \{2,3\}, \{2,4\}, \{2,5\}, \{3,4\}, \{3,5\}, \{4,5\}\}.
 $$
   Then we can easily check that there are at least 2 elements in $M_5$ which belong to at least half of the sets in $\cal{G}_2\cup \{G_{i_{2k+1}},\cdots, G_{i_m}\}$.

  \bigskip

  \quad\quad\quad\quad  (b.2.2) $n=3$. Now $\cal{G}_2=\{H_1,H_2,H_3\}$.

  \quad\quad\quad\quad\quad   (b.2.2.1) There exists $j\in \{1,2,3\}$ such that $H_j=\{1,2\}$. Then all the 2 elements in $\{1,2\}$ belong to at least 3 sets among 6 sets in $\cal{G}_2\cup \{G_{i_{2k+1}},\cdots, G_{i_m}\}$.
  \smallskip

\quad\quad\quad\quad\quad   (b.2.2.2) There exists $j\in \{1,2,3\}$ such that $H_j=\{1,3\}$. Then all the 2 elements in $\{1,3\}$ belong to at least 3 sets among 6 sets in $\cal{G}_2\cup \{G_{i_{2k+1}},\cdots, G_{i_m}\}$.
  \smallskip

\quad\quad\quad\quad\quad   (b.2.2.3) There exists $j\in \{1,2,3\}$ such that $H_j=\{1,4\}$. Then all the 2 elements in $\{1,4\}$ belong to at least 3 sets among 6 sets in $\cal{G}_2\cup \{G_{i_{2k+1}},\cdots, G_{i_m}\}$.
  \smallskip

    \quad\quad\quad\quad\quad    (b.2.2.4)  There exists $j\in \{1,2,3\}$ such that $H_j=\{4,5\}$. Without loss of generality, we assume that $H_1=\{4,5\}$ and thus $H_1\cup \{1,2,3\}=M_5$. Hence  it is enough to show that there exist two elements in $M_5$ which belong to at least half of the sets in $\{H_2,H_3\}\cup \{\{1,2,4\}, \{1,3,4\}\}$. Now in virtue of  (b.2.2.3), we need only consider
    \begin{eqnarray}\label{3.2}
    \{H_2,H_3\}\subset \{\{2,4\}, \{3,4\}, \{1,5\},\{2,5\}, \{3,5\}\},
    \end{eqnarray}
   If $\{2,4\}\in \{H_2,H_3\}$, then the 2 elements in $\{2,4\}$ belong to at least 2 sets among 4 sets in $\{H_2,H_3\}\cup \{\{1,2,4\}, \{1,3,4\}\}$.  If $\{3,4\}\in \{H_2,H_3\}$, then the 2 elements in $\{3,4\}$ belong to at least 2 sets among 4 sets in $\{H_2,H_3\}\cup \{\{1,2,4\},\{1,3,4\}\}$. Hence we need only consider the following 3 cases: $\{H_2,H_3\}=\{\{1,5\}, \{2,5\}\}, \{\{1,5\}, \{3,5\}\}, \{\{2,5\}, \{3,5\}\}\}$.  We can easily check that in all these 3 cases, there exist 2 elements in $M_5$ which belong to at least half of the sets in $\{H_2,H_3\}\cup \{\{1,2,4\},\{1,3,4\}\}$.
  \smallskip

   \quad\quad\quad\quad   (b.2.2.5)  There exists $j\in \{1,2,3\}$ such that $H_j=\{3,5\}$ or $H_j=\{2,5\}$. Similar to the analysis in (b.2.2.4), we can get that there exist two elements in $M_5$ which belong to at least half of the sets in $\cal{G}_2\cup \{G_{i_{2k+1}},\cdots, G_{i_m}\}$.

  \smallskip

   \quad\quad\quad\quad   (b.2.2.6)   $\cal{G}_2\subset \{\{1,5\}, \{2,3\}, \{2,4\},\{3,4\}\}$.  Now $ \cal{G}_2=\{\{2,3\}, \{2,4\}, \linebreak  \{3,4\}\}$.   We  can easily check that there exist two elements in $M_5$ which belong to at least half of the sets in $\cal{G}_2\cup \{G_{i_{2k+1}},\cdots, G_{i_m}\}$.

\bigskip

\quad\quad\quad\quad  (b.2.3) $n=4$.  Now $\cal{G}_2=\{H_1,H_2, H_3,H_4\}$ and  we have only the following 5 subcases:

\quad\quad\quad\quad\quad   (b.2.3.1) $\cal{G}_2=\{\{1,2\},\{1,3\}, \{1,4\},\{1,5\}\}$.

\quad\quad\quad\quad\quad   (b.2.3.2) $\cal{G}_2=\{\{2,1\},\{2,3\}, \{2,4\},\{2,5\}\}$.

\quad\quad\quad\quad\quad   (b.2.3.3) $\cal{G}_2=\{\{3,1\},\{3,2\}, \{3,4\},\{3,5\}\}$.

\quad\quad\quad \quad\quad   (b.2.3.4) $\cal{G}_2=\{\{4,1\},\{4,2\}, \{4,3\},\{4,5\}\}$.

\quad\quad\quad \quad\quad   (b.2.3.5) $\cal{G}_2=\{\{5,1\},\{5,2\}, \{5,3\},\{5,4\}\}$.

\smallskip

\quad As to (b.2.3.1), we know that $\{\{1,2,3,4\}, \{1,2,3,5\}, \{1,2,4,5\}, \{1,3,4,5\}\}\subset \cal{G}_4$. If $\cal{G}_4=\{\{1,2,3,4\}, \{1,2,3,5\}, \{1,2,4,5\}, \{1,3,4,5\}\}$, we can check that all the 4 elements in $\{1,2,3,4\}$ belong to at least half of the sets in $\cal{G}_2\cup \{G_{i_{2k+1}},\cdots, G_{i_m}\}\cup \cal{G}_4$ and thus belong to at least half of the sets in $\cal{F}$.
\smallskip

\quad As to (b.2.3.2), we can check that all the 2 elements in $\{1,2\}$ belong to at least half of the sets in $\cal{G}_2\cup \{G_{i_{2k+1}},\cdots, G_{i_m}\}$.

\smallskip

\quad As to (b.2.3.3), we can check that all the 2 elements in $\{1,3\}$ belong to at least half of the sets in $\cal{G}_2\cup \{G_{i_{2k+1}},\cdots, G_{i_m}\}$.

\smallskip

\quad As to (b.2.3.4), we can check that all the 2 elements in $\{1,4\}$ belong to at least half of the sets in $\cal{G}_2\cup \{G_{i_{2k+1}},\cdots, G_{i_m}\}$.

\smallskip

\quad As to (b.2.3.5), we can check that all the 2 elements in $\{1,5\}$ belong to at least half of the sets in $\cal{G}_2\cup \{G_{i_{2k+1}},\cdots, G_{i_m}\}$.

\bigskip

\quad (2.2.3.2-2),(2.2.3.2-4),(2.2.3.2-7),(2.2.3.2-8),(2.2.3.2-9): The proofs are similar to the one in  (2.2.3.2-1).

\bigskip

\quad (2.2.3.2-3):  We have two subcases: $|\cal{G}_2|=1$ and $|\cal{G}_2|\geq 2$.

\quad\quad (a) $|\cal{G}_2|=1$. Denote $\cal{G}_2=\{H\}$. Now all the 2 elements in $\{1,2\}$ belong to at least 3 sets among the four sets in $\cal{G}_2\cup \{G_{i_{2k+1}},\cdots, G_{i_m}\}$.

\quad\quad (b) $|\cal{G}_2|=n\geq 2$.  Denote $\cal{G}_2=\{H_1,\cdots,H_n\}$.  We have the following two subcases:
     \smallskip

  \quad\quad\quad  (b.1)  There exists a 2-element subset of $\{1,\cdots,n\}$ such that $H_i\cup H_j\in \cal{G}_4$, which implies that $H_i\cap H_j=\emptyset$. Then by Lemma \ref{lem-1.2}, we get that there are at least two elements in $M_5$ which belong to at least half of the sets in $\cal{F}$.
       \smallskip

        \quad\quad\quad  (b.2) For any $i,j=1,\cdots,n,i\neq j$,  $H_i\cup H_j\in \cal{G}_3$.  Now we can check that $n\leq 4$.  And so we have only three cases, which are $n=2,n=3,n=4$.

\quad\quad\quad\quad  (b.2.1) $n=2$.  Now all the 2 elements in $\{1,2\}$ belong to at least 3 sets among the 5 sets in $\cal{G}_2\cup \{G_{i_{2k+1}},\cdots, G_{i_m}\}$.

  \bigskip

  \quad\quad\quad\quad  (b.2.2) $n=3$. Now all the 2 elements in $\{1,2\}$ belong to at least 3 sets among the 6 sets in $\cal{G}_2\cup \{G_{i_{2k+1}},\cdots, G_{i_m}\}$.

\bigskip

\quad\quad\quad\quad  (b.2.3) $n=4$.  Now $\cal{G}_2=\{H_1,H_2, H_3,H_4\}$ and  we have only the following 5 subcases:

\quad\quad\quad\quad\quad   (b.2.3.1) $\cal{G}_2=\{\{1,2\},\{1,3\}, \{1,4\},\{1,5\}\}$.

\quad\quad\quad \quad\quad   (b.2.3.2) $\cal{G}_2=\{\{2,1\},\{2,3\}, \{2,4\},\{2,5\}\}$.

\quad\quad\quad\quad\quad   (b.2.3.3) $\cal{G}_2=\{\{3,1\},\{3,2\}, \{3,4\},\{3,5\}\}$.

\quad\quad\quad\quad\quad   (b.2.3.4)  $\cal{G}_2=\{\{4,1\},\{4,2\}, \{4,3\},\{4,5\}\}$.

\quad\quad\quad\quad\quad   (b.2.3.5) $\cal{G}_2=\{\{5,1\},\{5,2\}, \{5,3\},\{5,4\}\}$.

\smallskip

\quad As to (b.2.3.1) and (b.2.3.2) , we know that all the 2 elements in $\{1,2\}$ belong to at least 4 sets among the 7 sets in $\cal{G}_2\cup \{G_{i_{2k+1}},\cdots, G_{i_m}\}$.
\smallskip

\quad As to (b.2.3.3), we can check that all the 2 elements in $\{1,3\}$ belong to at least 4 sets among the 7 sets in $\cal{G}_2\cup \{G_{i_{2k+1}},\cdots, G_{i_m}\}$.

\smallskip

\quad As to (b.2.3.4), we can check that all the 2 elements in $\{1,4\}$ belong to at least 4 sets among the 7 sets in $\cal{G}_2\cup \{G_{i_{2k+1}},\cdots, G_{i_m}\}$.

\smallskip

\quad As to (b.2.3.5), we can check that all the 2 elements in $\{1,5\}$ belong to at least 4 sets among the 7 sets in $\cal{G}_2\cup \{G_{i_{2k+1}},\cdots, G_{i_m}\}$.

\bigskip

\quad (2.2.3.2-5) and  (2.2.3.2-6): The proofs are similar to the one in  (2.2.3.2-3).

\bigskip

\item[(2.2.3.3)]  $m-2k=4$. Now either $\{G_{i_{2k+1}},\cdots, G_{i_m}\}=\{\{1,2,3\}, \{1,2,4\},\{1,3,4\},\{2,3,4\}$ or $\{G_{i_{2k+1}},\cdots, G_{i_m}\}=\{\{1,2,3\}, \{1,2,5\},\{1,3,5\},\{2,3,5\}$. We only consider the former case. The proof for the latter case is similar.      We have two subcases: $|\cal{G}_2|=1$ and $|\cal{G}_2|\geq 2$.

\quad\quad (a) $|\cal{G}_2|=1$. Denote $\cal{G}_2=\{H\}$. Now all the 4 elements in $\{1,2,3,4\}$ belong to at least 3 sets among the 5 sets in $\cal{G}_2\cup \{G_{i_{2k+1}},\cdots, G_{i_m}\}$.

\quad\quad (b) $|\cal{G}_2|=n\geq 2$.  Denote $\cal{G}_2=\{H_1,\cdots,H_n\}$.  We have the following two subcases:
     \smallskip

  \quad\quad\quad  (b.1)  There exists a 2-element subset of $\{1,\cdots,n\}$ such that $H_i\cup H_j\in \cal{G}_4$, which implies that $H_i\cap H_j=\emptyset$. Then by Lemma \ref{lem-1.2}, we get that there are at least two elements in $M_5$ which belong to at least half of the sets in $\cal{F}$.
       \smallskip

        \quad\quad\quad  (b.2) For any $i,j=1,\cdots,n,i\neq j$,  $H_i\cup H_j\in \cal{G}_3$.  Now we can check that $n\leq 4$.  And so we have only three cases, which are $n=2,n=3,n=4$.

\quad\quad\quad\quad  (b.2.1) $n=2$.  Now all the 4 elements in $\{1,2,3,4\}$ belong to at least 3 sets among the 6 sets in $\cal{G}_2\cup \{G_{i_{2k+1}},\cdots, G_{i_m}\}$.

  \bigskip

  \quad\quad\quad\quad  (b.2.2) $n=3$.  Now $\cal{G}_2=\{H_1,H_2,H_3\}$. Note that all the 4 elements in $\{1,2,3,4\}$ belong to just 3 sets among the 4 sets in $\{\{1,2,3\}, \{1,2,4\},\{1,3,4\},\linebreak \{2,3,4\}$. Then we know that  if some element $j\in \{1,2,3,4\}$ belongs to at least one set in $\cal{G}_2$, then the element $j$ belong to at least 4 sets among 7 sets in $\cal{G}_2\cup \{\{1,2,3\}, \{1,2,4\},\{1,3,4\},\{2,3,4\}\}$.  By $|\cal{G}_2|=3$, we can easily know that  there exist at least 2 elements in $M_5$ which belong to at least half of the sets in $\cal{G}_2\cup \{G_{i_{2k+1}},\cdots, G_{i_m}\}$.

\bigskip

\quad\quad\quad\quad  (b.2.3) $n=4$.  Now $\cal{G}_2=\{H_1,H_2, H_3,H_4\}$ and  we have only the following 5 subcases:

\quad\quad\quad\quad\quad   (b.2.3.1) $\cal{G}_2=\{\{1,2\},\{1,3\}, \{1,4\},\{1,5\}\}$.

\quad\quad\quad \quad\quad   (b.2.3.2) $\cal{G}_2=\{\{2,1\},\{2,3\}, \{2,4\},\{2,5\}\}$.

\quad\quad\quad\quad\quad   (b.2.3.3) $\cal{G}_2=\{\{3,1\},\{3,2\}, \{3,4\},\{3,5\}\}$.

\quad\quad\quad\quad\quad   (b.2.3.4)  $\cal{G}_2=\{\{4,1\},\{4,2\}, \{4,3\},\{4,5\}\}$.

\quad\quad\quad\quad\quad   (b.2.3.5) $\cal{G}_2=\{\{5,1\},\{5,2\}, \{5,3\},\{5,4\}\}$.

\smallskip

\quad As to (b.2.3.1) - (b.2.3.4), we know that all the 4 elements in $\{1,2,3,4\}$ belong to at least 4 sets among the 8 sets in $\cal{G}_2\cup \{G_{i_{2k+1}},\cdots, G_{i_m}\}$.
\smallskip

\quad As to (b.2.3.5), we can check that all the 5 elements in $M_5$ belong to at least 4 sets among the 8 sets in $\cal{G}_2\cup \{G_{i_{2k+1}},\cdots, G_{i_m}\}$.
\end{itemize}
\end{itemize}
\end{itemize}
\end{itemize}
The whole proof is complete. \hfill\fbox

\section{A counterexample for $S_1$-Frankl conjecture, some questions  and  remarks}

\subsection{A counterexample and three questions}

The following example shows that $S_1$-Frankl conjecture claims too much and is not valid in general, which is introduced by an anonymous referee.

Let $n=9,M_9=\{1,2,\ldots,9\}$, and
\begin{eqnarray*}
\cal{F}&=&\{\emptyset, \{1,2,7,8\}, \{3,4,7,9\}, \{5,6,8,9\},\\
&&\{1,2,7,8,9\}, \{3,4,7,8,9\}, \{5,6,7,8,9\},\\
&&\{1,2,3,4,7,8,9\}, \{1,2,5,6,7,8,9\}, \{3,4,5,6,7,8,9\},\\
&&\{1,2,3,4,5,6,7,8,9\}\}.
\end{eqnarray*}
The family $\cal{F}$ is union-closed, $|\cal{F}|=11,T(\cal{F})=4$. Each of the elements in $\{1,2,3,4,5,6\}$ appears 5 sets in $\cal{F}$, and only 3 elements 7, 8 and 9 appear in more than half of the sets in $\cal{F}$.

Based on the above example, we can easily construct example in which $T(\cal{F})\geq 5$, and the result in $S_1$-Frankl conjecture is not true. Three  natural questions are as follows:

{\it Question 1. If $T(\cal{F})=3$, does there exist 3 elements which belong to at least half of the sets in $\cal{F}$?}

{\it Question 2. If $T(\cal{F})=2$, does there exist 2 elements which belong to at least half of the sets in $\cal{F}$?}

{\it Question 3. If $T(\cal{F})>\frac{n}{2}$, does there exist $T(\cal{F})$ elements which belong to at least half of the sets in $\cal{F}$, i.e. does the claim in $S_1$-Frankl conjecture hold in this case?}

Obviously, Question 2 is a part of $S_2$-Frankl conjecture.

\subsection{Some remarks}

Although $S_1$-Frankl conjecture is not valid in general, we still believe that the idea of using the structure of $\cal{F}$ may be used to study the original Frankl's conjecture.

Up to now, we do not know whether Frankl's conjecture holds when $n\geq 13$. In the following, take $n=14$ for example. Let $M_{14}=\{1,2,\cdots,14\}$ and $\cal{F}\subset 2^{M_{14}}=\{A: A\subset M_{14}\}$ with $\cup_{A\in\cal{F}}A=M_{14}$.  Suppose that $\cal{F}$ is union-closed and  $\emptyset \in \cal{F}$.  For any $k=1,2,\cdots,14$, define
$$
\cal{M}_k=\{A\in 2^{M_{14}}: |A|=k\},\ \ n_k=|\cal{F}\cap \cal{M}_k|,
$$
and
$$
T(\cal{F})=\inf\{1\leq k\leq 14: n_k>0\}.
$$
Then $1\leq T(\cal{F})\leq 14$.  If $T(\cal{F})=14$, then $\cal{F}=\{\emptyset, M_{14}\}$.  Now all the 14 elements in $M_{14}$ belongs to half of the sets in $\cal{F}$.

If $T(\cal{F})\in \{7,\cdots,13\}$, we can easily prove that there exists an element in $M_{14}$ which belongs to at least half of the sets in $\cal{F}$. Without loss of generality, we assume that $T(\cal{F})=7$. By
$$
\frac{1}{14}\sum_{k=7}^{13}k\times n_k\geq \frac{7(n_7+\cdots+n_{13})}{14}= \frac{n_7+\cdots+n_{13}}{2},
$$
we know that there exists at least an element in $M_{14}$ which belongs to at least half of the sets in $\cal{F}\cap \left(\cup_{k=7}^{13}\cal{M}_k\right) $ and thus belongs to at least half of the sets in $\cal{F}$.

Hence we need only to consider the four cases: $T(\cal{F})=3,4,5,6$.  In the following, we give some (partial, not complete) discussions for the case $T(\cal{F})=6$ for the illustration.

\begin{itemize}
\item[(1)] $T(\cal{F}\backslash\cal{M}_6)=14$.  Now $\cal{F}=\{\emptyset,M_{14}\}\cup \cal{G}_6$, where $\cal{G}_6$ is a nonemptyset subset of $\cal{M}_6$.

     For any two different sets $A$ and $B$ in $\cal{M}_6$, we have that $7\leq |A\cup B|\leq 12$. Hence $n_6=1$. Denote $\cal{G}_6=\{G\}$. Then all the 6 elements in $G$ belong to at least half of the sets in $\cal{F}$.

\item[(2)] $T(\cal{F}\backslash\cal{M}_6)=13$. Now $\cal{F}=\{\emptyset,M_{14}\}\cup \cal{G}_6\cup \cal{G}_{13}$, where $\cal{G}_i$ is a nonemptyset subset of $\cal{M}_i$ for
$i\in \{6, 13\}$, and thus $n_6\geq 1,n_{13}\geq 1$.

By the analysis in (1), we know that $n_6=1$. By
$$
\frac{6+13n_{13}}{14}=\frac{(6+6n_{13})+7n_{13}}{14}\geq \frac{12+7n_{13}}{14} >\frac{n_6+n_{13}}{2},
$$
 we get that there exists at least one element in $M_{14}$ which belongs to at least half of the sets in $\cal{G}_6\cup \cal{G}_{13}$ and thus belongs to at least half of the sets in $\cal{F}$.

\item[(3)] $T(\cal{F}\backslash\cal{M}_6)=12$. Now we have two subcases:

\begin{itemize}
\item[(3.1)]  $T(\cal{F}\backslash(\cal{M}_6\cup \cal{M}_{12}))=14$.  Now $\cal{F}=\{\emptyset,M_{14}\}\cup \cal{G}_6\cup \cal{G}_{12}$.

\begin{itemize}
\item[(3.1.1)] $n_6=1$. By
$$
\frac{6+12n_{12}}{14}=\frac{(6+5n_{12})+7n_{12}}{14}\geq \frac{11+7n_{12}}{14}
>\frac{n_6+n_{12}}{2},
$$
we get that there exists at least one element in $M_{14}$ which belongs to at least half of the sets in $\cal{G}_6\cup \cal{G}_{12}$ and thus belongs to at least half of the sets in $\cal{F}$.

\item[(3.1.2)] $n_6\geq 2$.  Denote $\cal{G}_6=\{H_1,\cdots,H_{n_6}\}$. Then, for any $i,j=1,\cdots,n_6,i\neq j$, we must have $|H_i\cup H_j|=12$, which implies that $H_i\cap H_j=\emptyset$. It follows that $n_6=2$.  By
    $$
    \frac{6\times 2+12 n_{12}}{14}=\frac{(12+5n_{12})+7 n_{12}}{14}\geq \frac{17+7 n_{12}}{14}>\frac{n_6+n_{12}}{2},
        $$
   we get that there exists at least one element in $M_{14}$ which belongs to at least half of the sets in $\cal{G}_6\cup \cal{G}_{12}$ and thus belongs to at least half of the sets in $\cal{F}$.
\end{itemize}

\item[(3.2)]  $T(\cal{F}\backslash(\cal{M}_6\cup \cal{M}_{12}))=13$. Now $\cal{F}=\{\emptyset,M_{14}\}\cup \cal{G}_6\cup \cal{G}_{12}\cup \cal{G}_{13}$.

For any two different sets $A,B$  in $\cal{M}_6$, we have that $7\leq |A\cup B|\leq 12$. Then by following the proof in (3.1), we can get that there exists an element in  $M_{14}$ which belongs to at least half of the sets in $\cal{F}$. In fact, we have the following two subcases:

\begin{itemize}
\item[(3.2.1)] $n_6=1$. Now we have
$$
\frac{6+12n_{12}+13n_{13}}{14}=\frac{(6+5n_{12}+6n_{13})+7n_{12}+7n_{13}}{14}>\frac{n_6+n_{12}+n_{13}}{2}.
$$

\item[(3.2.2)] $n_6=2$. Now we  have
$$
\frac{6\times 2+12n_{12}+13n_{13}}{14}=\frac{(12+5n_{12}+6n_{13})+7n_{12}+7n_{13}}{14}>\frac{n_6+n_{12}+n_{13}}{2}.
$$
\end{itemize}
\end{itemize}
\end{itemize}

%
%
%

\vskip 0.5cm
{ \noindent {\bf\large Acknowledgments}

\noindent We thank an anonymous referee very much  for giving a counterexample for $S_1$-Frankl conjecture to the previous version of this paper. This work was supported by National Natural Science Foundation of China (Grant No. 11771309) and the Fundamental Research Funds for the Central Universities.


\begin{thebibliography}{12}

\bibitem{BM08} Bo\v{s}njak I., Markovi\'{c} P.:  Then 11-element case of Frankl's conjecture, The Electronic J. Comb. {\bf 15}, \#R88 (2008)

\bibitem{BS15} Bruhn H., Schaudt O.: The journey of the union-closed sets conjecture, Graphs and Comb. {\bf 31},  2043-2074 (2015)




\bibitem{GY98} Gao W., Yu H.: Note on the union-closed sets conjecture, Ars Combin. {\bf 49}, 280-288 (1998)

\bibitem{JV98} Johnson R. T.,  Vaughan T. P.: On union-closed families, I. J. Combin. Theory Ser. A {\bf  84(2)}, 242-249 (1998)

\bibitem{Lo94-a} Lo Faro G.: A note on the union-closed sets conjecture, J. Austral. Math. Soc. Ser A {\bf 57(2)}, 230-236 (1994)

\bibitem{Lo94-b}  Lo Faro G.: Union-closed sets conjecture: improved bounds, J. Comb. Math. Comb. Comput. {\bf 16}, 97-102  (1994)

\bibitem{Ma07} Markovi\'{c} P.: An attempt at Frankl's conjecture, Proceedings of the Novi Sad Algebraic Conference 2005 (NSAC'05), a special issue of Publ. Math. Inst. (Beograd) (N. S.) {\bf 81(95)}, 29-43 (2007)

\bibitem{Mo06} Morris R.: FC-families and improved bound for Frankl's conjecture, European J. Combin. {\bf 27(2)}, 269-282 (2006)

\bibitem{Po92} Poonen B.: Union-closed families, J. Comb. Theory, Ser. A {\bf 59(2)}, 253-268 (1992)

\bibitem{Ri85} Rival I. (Ed.): Graphs and Order,  Reidel, Dordrecht/Boston (1985)

\bibitem{RS10} Roverts I., Simpson J.: A note on the union-closed sets conjecture, Austral. J. Comb. {\bf 47}, 265-267 (2010)


\bibitem{SR89}  Sarvate  D.G.,  Renaud J.-C.: On the union-closed sets conjecture, Ars Combin. {\bf 27}, 149-154 (1989)

\bibitem{St86}  Stanley R. P., Enumerative Combinatorics, Vol. I, Wadsworth \& Brooks/Cole, Belmont, CA (1986)

\bibitem{Va02} Vaughan T. P.: Families implying the Frankl conjecture, European J. Combin. {\bf 23(7)}, 851-860 (2002)

\bibitem{Va03} Vaughan T. P.: A note on the union-closed sets conjecture, J. Combin. Math. Combin. Comput. {\bf 45}, 95-108 (2003)

\bibitem{Va04} Vaughan T. P.: Three-sets in a union-closed family, J. Combin. Math. Combin. Comput. {\bf 49}, 73-84 (2004)


\bibitem{ZV12} Zivkovi\'{c} M.,  Vu\v{c}kovi\'{c} B.:  The 12 element case of Frankl's conjecture, Preprint (2012)



\end{thebibliography}
\end{document}